\documentclass[11pt]{article}
\usepackage{myarticle-macros,graphicx}
\usepackage{subfig}

\newcommand{\cplus}{\mathbb{C}^+}
\newcommand{\rplus}{\mathbb{R}^+}
\newcommand{\mapa}{Mar\v{c}enko-Pastur\ }

\newcommand{\tr}[1]{\text{trace}\left(#1\right)}
\newcommand{\myreal}[1]{\text{Re}\left(#1\right)}
\newcommand{\imag}[1]{\text{Im}\left(#1\right)}

\newtheorem*{fact*}{Fact}
\newcommand{\WassDist}[2]{d_W\left(#1,#2\right)}

\newcommand{\vinfty}{v_{\infty}}

\newcommand{\Hp}{H_{p}}
\newcommand{\hatHp}{\widehat{H}_p}
\newcommand{\Hinfty}{H_{\infty}}
\newcommand{\hatHinfty}{\widehat{H}_{\infty}}
\newcommand{\discHinfty}{H_{M_n}}
\newcommand{\vmin}{u_{\mathrm{min}}}
\newcommand{\cplusabovevmin}{\cplus \cap \left\{\imag{z}>\vmin\right\}}
\newcommand{\Fp}{F_{p}}
\newcommand{\tildeFp}{\widetilde{F}_p}
\newcommand{\Finfty}{F_{\infty}}
\newcommand{\esd}{\Fp}
\newcommand{\stesd}{m_{\esd}}
\newcommand{\stesdd}{v_{\esd}}
\newcommand{\EmpCovMat}{S_{p}}

\DeclareMathOperator{\argmin}{argmin}
\title{Spectrum estimation for large dimensional covariance matrices
using random matrix theory}
\usepackage{geometry}
\usepackage{setspace}
\author{Noureddine El Karoui \thanks{\textbf{Acknowledgements:} The author is grateful to Alexandre d'Aspremont,
Peter Bickel, Laurent El Ghaoui, Elizabeth Purdom, John Rice,
Saharon Rosset and Bin Yu for stimulating discussions and comments
at various stages of this project. Support from NSF grant
DMS-0605169 is gratefully acknowledged.
\textbf{AMS 2000 SC: } Primary 62H12, Secondary 62-09.
\textbf{Key words and Phrases~: } covariance
matrices, principal component analysis, eigenvalues of covariance
matrices, high-dimensional inference, random matrix theory,
Stieltjes transforms, \mapa equation, convex optimization.
 \textbf{Contact~:} \texttt{nkaroui@stat.berkeley.edu}} \\
\textit{Department of Statistics,}\\ \textit{University of California, Berkeley} }
\begin{document}
\maketitle
\begin{abstract}
Estimating the eigenvalues of a population covariance matrix from a
sample covariance matrix is a problem of fundamental importance in
multivariate statistics; the eigenvalues of covariance matrices play
a key role in many widely techniques, in particular in Principal
Component Analysis (PCA). In many modern data analysis problems,
statisticians are faced with large datasets where the sample size,
$n$, is of the same order of magnitude as the number of variables
$p$. Random matrix theory predicts that in this context, the
eigenvalues of the sample covariance matrix are not good estimators
of the eigenvalues of the population covariance.

We propose to use a fundamental result in random matrix theory, the
\mapa equation, to better estimate the eigenvalues of large
dimensional covariance matrices. The \mapa equation holds in very
wide generality and under weak assumptions. The estimator we obtain
can be thought of as ``shrinking" in a non linear fashion the
eigenvalues of the sample covariance matrix to estimate the
population eigenvalue. Inspired by ideas of random matrix theory, we
also suggest a change of point of view when thinking about
estimation of high-dimensional vectors: we do not try to estimate
directly the vectors but rather a probability measure that describes
them. We think this is a theoretically more fruitful way to think
about these problems.

Our estimator gives fast and good or very good results in extended
simulations. Our algorithmic approach is based on convex
optimization. We also show that the proposed estimator is
consistent.
\end{abstract}
\section{Introduction}
With data acquisition and storage now easy, today's statisticians
often encounter datasets for which the sample size, $n$ and the
number of variables $p$, are both large: in the order of hundreds,
thousands, millions, or even billions in situations such as web
search problems.

The analysis of these datasets using classical methods of
multivariate statistical analysis requires some care. While the
ideas are still relevant, the intuition for the estimators that are
used and the interpretation of the results are often - implicitly -
justified by assuming an asymptotic framework of $p$ fixed and $n$
growing infinitely large. This assumption was consistent with the
practice of statistics when these ideas were developed, since
investigation of datasets with a large number of variables was very
difficult. A better theoretical framework for modern - i.e large $p$
- datasets, however is the assumption of the so-called ``large $n$,
large $p$" asymptotics. In other words, one should consider that
both $n$ and $p$ go to infinity, perhaps with the restriction that
their ratio goes to a finite limit $\gamma$, and draw practical
insights from the theoretical results obtained in this setting.

We will turn our attention to an object of central interest in
multivariate statistics: the eigenvalues of covariance matrices. A
key application is  Principal Components Analysis (PCA), where one
searches for a good low dimensional approximation to the data by
projecting the data on the ``best" possible $k$ dimensional
subspace: here ``best" means that the projected data explain as much
variance in the original data as possible. This amount of variance
explained is measured by the eigenvalues of the population
covariance matrix, $\Sigma_p$, and hence we need to find a way to
estimate those eigenvalues. We will discuss in the course of the
paper other problems where the eigenvalues of $\Sigma_p$ play a key
role.

We take a moment here to give a few examples that illustrate the
differences that occur under the different asymptotic settings. To
pose the problem more formally, let us say that we observe iid
random vectors $X_1,\ldots,X_n$ in $\mathbb{R}^p$, and that the
covariance of $X_i$ is $\Sigma_p$. We call $X$ the data matrix whose
rows are the $X_i$'s. In the classical context, where $p$ is fixed
and $n$ goes to $\infty$, a fundamental result of
\citep{anderson63} says that the eigenvalues of the sample covariance
matrix $S_p=(X-\bar{X})'(X-\bar{X})/(n-1)$ are good estimators of
the population eigenvalues (i.e the eigenvalues of $\Sigma_p$). More
precisely, calling $l_i$ the ordered eigenvalues of $S_p$ ($l_1\geq
l_2\ldots$) and $\lambda_i$ the ordered eigenvalues of $\Sigma_p$
($\lambda_1\geq \lambda_2 \ldots$), it was shown in
\citep{anderson63} that
$$
\sqrt{n}(l_i-\lambda_i)\Rightarrow {\cal N}(0,2\lambda_i^2)\;,
$$
when the $X_i$ are normally distributed and all the $\lambda_i$'s
are distinct. This result provided rigorous grounds for estimating
the eigenvalues of the population covariance matrix, $\Sigma_p$,
with the eigenvalues of the sample covariance matrix, $S_p$, when
$p$ is small compared to $n$. (For more details on Anderson's
theorem, we refer the reader to
\citep{anderson03} Theorem 13.5.1.)

Shifting assumptions to ``large $n$, large $p$" asymptotics induces
fundamental differences in the behavior of multivariate statistics,
some of which we will highlight in the course of the paper. As a
first example, let us consider the case where $\Sigma_p=\id_p$, so
all the population eigenvalues are equal to $1$. A result first
shown in
\citep{geman80} under some moment growth assumptions, and later refined in
\citep{yinbaikrishnaiah88}, states that if the entries of the $X_i$'s are i.i.d and  have
a fourth moment, and if $p/n\tendsto \gamma$, then
$$
l_1\tendsto (1+\sqrt{\gamma})^2 \;\; \text{a.s.}
$$
In particular, $l_1$ is not a consistent estimator of $\lambda_1$.
Note that by picking $n=p$, $l_1$ tends to 4 whereas $\lambda_1= 1$.
(For more general $\Sigma_p$, see \citep{nekGencov} Section 4.3 for
numerically explicit results about the limit of $l_1$.)

As the case of $\Sigma_p=\id_p$ illustrated, when $n$ and $p$ are
both large, the largest sample eigenvalue is biased, sometime
dramatically so. Hence, we should correct this bias in the largest
sample eigenvalue(s) if we want to use them in data analysis.
Theoretical results predict that the behavior of extreme sample
eigenvalues can be quite subtle; in particular, depending on how far
an isolated population eigenvalue is from the bulk of the population
spectrum, the corresponding sample eigenvalue can either be
isolated, and far away from the bulk of the sample eigenvalues, or
be absorbed by the bulk of the sample eigenvalues (see
\citep{bbap}, \citep{nekGencov}, \citep{baiksilverstein04}, \citep{debashis}).
One thing is however clear from the most recent theoretical results
: if we wish to de-bias extreme sample eigenvalues, we need an
accurate estimate of the so-called population spectral distribution,
a probability measure that characterizes the population eigenvalues
(see \citep{nekGencov}). This is what our algorithm will deliver.

We have so far mostly discussed extreme sample eigenvalues. However,
much is also known about the behavior of the whole vector of sample
eigenvalues $(l_1,l_2,\ldots,l_p)$ and its asymptotic behavior. In
particular, theory predicts that in the ``large $n$, large $p$"
case, the scree plot (i.e the plot of the sample eigenvalues vs.
their rank; see
\citep{mardiakentbibby}) becomes uninformative and deceptive. What we
propose in this paper is to use random matrix theory to develop
practically useful tools to remedy the flaws appearing in some
widely used tools in multivariate statistics.

Before we discuss how we will go about it, let us briefly discuss
some issues that arise when estimating vectors of large dimension,
since working in an asymptotic setting where $p\tendsto \infty$ is
not without additional difficulties. Since we will try to estimate
vectors of increasingly larger and larger size, an appropriate
notion of convergence is needed if we want to quantify the quality
of our estimators. Standard norms in high-dimensions not necessarily
a very good choice: for instance, if we are in $\mathbb{R}^{100}$,
and make an error of size 1/100 in all coordinates, the resulting
$l_1$ error is 1, even though, at least intuitively, it would seem
like we are doing well. Also, if we made a large error (say size 1)
in one direction, the $l_2$ norm would be large (larger than 1 at
least), even though we may have gotten the structural information
about this vector (and almost all its coordinates) ``right".
Inspired by ideas of random matrix theory, we propose to associate
to high-dimensional vectors probability measures that describe them.
We will  explain this in more detail in Section
\ref{subsec:changePointView}.
After this change of point of view, our focus becomes trying to
estimate these measures. Why choosing to estimate measures? The
reasons are many. Chief among them is that this approach will allow
us to look into the structure of the population eigenvalues. For
instance, we would like to be able to say whether all population
eigenvalues are equal, or whether they are clustered around say two
values, or if they are uniformly spread out on an interval. Because
the ratio $p/n$ can make the scree plot appear smooth (and hence in
some sense uninformative) regardless of the true population
eigenvalue structure, this structural information is not well
estimated by currently existing methods. We discuss other practical
benefits (like scalability with $p$) of the measure estimation
approach in
\ref{subsubsec:ChooseEstimMeas}. In the context of PCA, where
usually the concern is not to estimate each population eigenvalues
with very high precision, but rather to have an idea of the
structure of the population spectrum to guide the choice of
lower-dimensional subspaces on which to project the data, this
measure approach is particularly appealing. Examples to come later
in the paper will illustrate this point.

Random matrix theory plays a key role in our approach to this
measure estimation problem. A main ingredient of our method is a
fundamental result, which we call the
\mapa equation (see Theorem
\ref{thm:mapa}), which relates the asymptotic behavior of the sample
eigenvalues to the population eigenvalues. The assumptions under
which the theorem holds are very weak (a fourth moment condition)
and hence it is very widely applicable. Until now, this theorem has
not been used to do inference on population eigenvalues. Partly this
is because in its general form it has not received much attention in
statistics, and partly because the inverse problem that needs to be
considered is very hard to solve if it is not posed the right way.
We propose an original way to approach inverting the
\mapa equation. In particular, we will be able to estimate
given the eigenvalues of the sample covariance matrix $S_p$ the
probability measure, $\Hp$, that describes the population
eigenvalues. We use the standard names \textit{empirical spectral
distribution} for $\Fp$ and \textit{population spectral
distribution} for $\Hp$. It is important to state clearly what
asymptotic framework we place ourselves in. We will consider that
when $p$ and $n$ go to infinity, $\Hp$ stays fixed. In particular,
it has a limit, denoted $\Hinfty$. We call this framework
``asymptotics at fixed spectral distribution". Of course, fixing
$\Hp$ does not imply that we fix $p$. For instance, sometime we will
have $\Hp=\delta_1$, for all $p$. Since the parameter of interest in
our problems is really the measure $\Hp$, the fixed spectral
distribution asymptotics corresponds to classical assumptions for
parameter estimation in statistics, where the parameter does not
change with the number of variables observed. We refer the reader to
\ref{subsubsec:IsoEigen} for a more detailed discussion.

To solve the inverse problem posed by the \mapa equation, we propose
to discretize the
\mapa equation and then use convex optimization methods to solve the
discretized version of the problem. In doing so, we obtain a fast
and provably accurate algorithm to estimate the population parameter
of interest, $\Hp$, from the sample eigenvalues. The approach is
non-parametric since no assumptions are made \textit{a priori} on
the structure of the population eigenvalues. One outcome of the
algorithm is an efficient graphical method to look at the structure
of the population eigenvalues. Another outcome is that since we have
an estimate of the measure that describes the population
eigenvalues, standard statistical ideas then allow us to get
estimates of the individual population eigenvalues $\lambda_i$. Some
subtle problems may arise when doing so and we address them in
\ref{subsubsec:IsoEigen}. The final result of the algorithm can be
thought of as performing non-linear shrinkage of the sample
eigenvalues to estimate the population eigenvalues.

We want to highlight two contributions of our paper. First, we
propose to estimate measures associated with high-dimensional
vectors rather than estimating the vectors. This gives rise to
natural notions of consistency and accuracy of our estimates which
are reasonable theoretical requirements for any estimator to
achieve. And second, we make use, for the first time, of a
fundamental result of random matrix theory to solve an important
practical problem in multivariate statistics.

The rest of the paper is divided into four parts. In Section
\ref{sec:RMTbackground}, we give some background on results in
Random Matrix Theory that will be needed. We do not assume that the
reader has any familiarity with the topic.  In Section
\ref{sec:Algo}, we present our algorithm to estimate $\Hp$, the population
spectral distribution, and also the population eigenvalues. In
Section
\ref{sec:Simulations}, we present the results of some simulations.
We give in Section
\ref{sec:Consistency} a proof of consistency of our algorithm. The
Appendix contains some details on implementation of the algorithm.

A note on notation is needed before we start: in the rest of the
paper, $p$ will always be a function of $n$, with the property that
$p(n)/n\tendsto \gamma$ and $\gamma\in(0,\infty)$. To avoid
cumbersome notations, we will usually write $p$ and not $p(n)$.
\section{Background: Random matrix theory of sample covariance
matrices}\label{sec:RMTbackground} There is a large body of work
concerned with the limiting behavior of the eigenvalues of a sample
covariance matrix when $p$ and $n$ both go to $\infty$; it
constitutes an important subset of what is commonly known as Random
Matrix Theory, to which we now turn. This is a wide area of
research, of which we will only give a very quick and self-contained
overview. Our eventual aim in this section is to introduce a
fundamental result, the \mapa equation, that relates the asymptotic
behavior of the eigenvalues of the sample covariance matrix to that
of the population covariance in the ``large $n$, large $p$"
asymptotic setting. The formulation of the result requires that we
introduce some concepts and notations.
\subsection{Changing point of views: from vectors to measures}
\label{subsec:changePointView}
One of the first problems to tackle is to find a mathematically
efficient way to express the limit of a vector whose size grows to
$\infty$. (Recall that there are $p$ eigenvalues to estimate in our
problem and $p$ goes to $\infty$.) A fairly natural way to do so is
to associate to any vector a probability measure. More explicitly,
suppose we have a vector $(y_1,\ldots,y_p)$ in $\mathbb{R}^p$. We
can associate to it the following measure:
$$
dG_p(x)=\frac{1}{p}\sum_{i=1}^p \delta_{y_i}(x)\;.
$$
$G_p$ is thus a measure with $p$ point masses of equal weight, one
at each of the coordinates of the vector.

In the rest of the paper, we will denote by $H_p$ the spectral
distribution of the population covariance matrix $\Sigma_p$, i.e the
measure associated with the vector of eigenvalues of $\Sigma_p$. We
will refer to $H_p$ as \emph{the population spectral distribution}.
We can write this measure as
$$
dH_p(x)=\frac{1}{p}\sum_{i=1}^p \delta_{\lambda_i}(x)\;,
$$
where $\delta_{\lambda_i}$ is a point mass, of mass 1, at
$\lambda_i$. We also call $\delta_{\lambda_i}$ a ``dirac" at
$\lambda_i$. The simplest example of population spectral
distribution is found when $\Sigma_p=\id_p$. In this case, for all
$i$, $\lambda_i=1$, and $dH_p=\delta_1$. So the population spectral
distribution is a point mass at 1 when $\Sigma_p=\id_p$.

Similarly, we will denote by $\Fp$ the measure associated with the
eigenvalues of the sample covariance matrix $\EmpCovMat$. We refer
to $\Fp$ as the \emph{empirical spectral distribution}.
Equivalently, we define
$$
d\esd(x)=\frac{1}{p}\sum_{i=1}^p \delta_{l_i}(x)\;.
$$

The change of focus from vector to measure implies a change of focus
in the notion of convergence we will consider adequate. In
particular, for consistency issues, the notion of convergence we
will use is weak convergence of probability measures. While this is
the natural way to pose the problem mathematically, we may ask if it
will allow us to gather the statistical information we are looking
for. An example of the difficulties that arise is the following.
Suppose $dH_p=(1-1/p)
\, \delta_1 + 1/p \, \delta_2$. In other words, the population
covariance has one eigenvalue that is equal to 2 and $(p-1)$ that
are equal to 1. Clearly, when $p\tendsto \infty$, $H_p$ weakly
converges to $\Hinfty$, with $d\Hinfty=\delta_1$. So all information
about the large and isolated eigenvalue $2$, which is present in
$\Hp$ for all $p$ and is naturally of great interest in PCA, seems
lost in the limit. This is not the case when one does asymptotic at
fixed spectral distribution and consider that we are following a
sequence of models which are going to infinity with
$H_p=H_{p_0}=\Hinfty$, where $p_0$ is the $p$ which is given by the
data set. Fixed distribution asymptotics is more akin to what is
done in classical statistics and we place ourselves in this
framework. We refer the reader to \ref{subsubsec:IsoEigen} for a
more detailed justification of our point.

In other respects, associating a measure to a vector in the way we
described is meaningful mostly when one wants to have information
about the whole set of values taken by the coordinates of the
vector, and not about each coordinate. In particular, when going
from vector to measure as described above we are losing all
coordinate information: permuting the coordinates would drastically
change the vector but yield the same measure. However, in the case
of vectors of eigenvalues, since there is a canonical way to
represent the vector (the $i$-th largest eigenvalue occupying the
$i$-th coordinate), the information contained in the measure is
sufficient. This measure approach is especially good when we are not
focused on getting all the fine details of the vectors right, but
rather when we are looking for structural information concerning the
values taken by the coordinates.

An important area of random matrix theory for sample covariance
matrices is concerned with understanding the properties of $\esd$ as
$p$ (and $n$) go to $\infty$. A key theorem , which we review later
(see Theorem \ref{thm:mapa}), states that for a wide class of sample
covariance matrices, $\Finfty$, the limit of $\esd$, is
asymptotically non-random. Furthermore, the theorem connects
$\Finfty$ to $\Hinfty$, the limit of $H_p$: given $\Hinfty$, we can
theoretically compute $\Finfty$, by solving a complicated equation.
In data analysis, we observe the empirical spectral distribution,
$\Fp$. Our goal, of course, as far as eigenvalues are concerned, is
to estimate the population spectral distribution, $\Hp$. Our method
will ``invert" the relation between $\Finfty$ and $\Hinfty$, so that
we can go from $\Fp$ to $\hatHp$,  an estimate of $\Hp$. The method
does not work directly with $\Fp$ but with a tool that is similar in
flavor to the characteristic function of a distribution: the
Stieltjes transform of a measure. We introduce this tool in the next
subsection. As we will see later, it will also play a key role in
our algorithm.

\subsection{The Stieltjes transform of measures}
A large number of results concerning the asymptotic properties of
the eigenvalues of large dimensional random matrices are formulated
in terms of limiting behavior of the Stieltjes transform of their
empirical spectral distributions. The Stieltjes transform is a
convenient and very powerful tool in the study of the convergence of
spectral distribution of matrices (or operators), just as the
characteristic function of a probability distribution is a powerful
tool for central limit theorems. Most importantly, there is a simple
connection between the Stieltjes transform of the spectral
distribution of a matrix and its eigenvalues.

By definition, the Stieltjes transform of a measure $G$ on
$\mathbb{R}$ is defined as
$$
m_G(z)=\int \frac{dG(x)}{x-z}\;, \text{ for }z\in \cplus ,
$$
where $\cplus\triangleq
\mathbb{C}\bigcap\{z:\, \imag{z}>0\}$ is the set of complex numbers with strictly
positive imaginary part. The Stieltjes transform appears to be known
under several names in different areas of mathematics. It is
sometimes referred to as Cauchy or Abel-Stieltjes transform. Good
references about Stieltjes transforms include \citep[Sections
3.1-2]{akhiezer65},
\citep[Chapter 32]{lax}, \citep[Chapter 3]{hiaipetz00} and \citep{geronimohill03}.

For the purpose of this paper, where will consider only compactly
supported measures, the following results will be needed:
\begin{fact*}
Important properties of Stieltjes transforms of measures on
$\mathbb{R}$:
\begin{enumerate}
\item If $G$ is a probability measure, $m_G(z)\in \cplus$  if $z\in
\cplus$ and $\lim_{y\tendsto \infty} -iy m_G(iy) =1$.
\item If $F$ and $G$ are two measures, and if $m_F(z)=m_G(z)$, for all $z \in \cplus$,
then $G=F$, a.e.
\item \citep[Theorem 1]{geronimohill03}: If $G_n$ is a sequence of probability measures and
$m_{G_n}(z)$ has a (pointwise) limit $m(z)$ for all $z \in \cplus$,
then there exists a probability measure $G$ with Stieltjes transform
$m_G=m$ if and only if $\lim_{y\tendsto
\infty} -iy m(iy) =1$.  If it is the case,  $G_n$ converges weakly to $G$.
\item \citep[Theorem 2]{geronimohill03}: The same is true if the convergence happens only for an
infinite sequence $\{z_i\}_{i=1}^{\infty}$ in $\cplus$ with a limit
point in $\cplus$.
\item If $t$ is a continuity
point of the cdf of $G$, $dG(t)/dt=\lim_{\eps\tendsto 0}
\frac{1}{\pi}\imag{m_G(t+i \eps)}$
\end{enumerate}
For proofs, we refer the reader to \citep{geronimohill03}.
\end{fact*}

Note that the Stieltjes transform of the spectral distribution
$\Gamma_p$ of a $p\times p$ matrix $A_p$ is just
$$
m_{\Gamma_p}(z)=\frac{1}{p}\tr{(A_p-z\id_p)^{-1}}\;.
$$
Finally, it is clear that points 3 and 4 above can be used to show
convergence of probability measures if one can control the
corresponding Stieltjes transforms.

\subsection{A fundamental result: the \mapa equation}
In the study of covariance matrices, a remarkable result exists that
describes the limiting behavior of the empirical spectral
distribution, $\Finfty$, in terms of the limiting behavior of the
population spectral distribution, $\Hinfty$. The connection between
these two measures is made through an equation that links the
Stieltjes transform of the empirical spectral distribution to an
integral against the population spectral distribution. We call this
equation the \mapa equation because it first appeared in the
landmark paper of
\citep{mp67}. The result was independently re-discovered  in
\citep{wachter78} and then refined in
\citep{silversteinbai95} and
\citep{silverstein95}. In particular, \citep{silverstein95} is the
only paper where the case of a non-diagonal population covariance is
tackled.

In what follows, we will be working with an $n\times p$ data matrix
$X$. We call $\EmpCovMat=X^*X/n$ and denote $\stesd$ the Stieltjes
transform of the spectral distribution, $\esd$, of $\EmpCovMat$. We
will call $v_{\esd}$ the function defined by
$v_{\esd}(z)=(1-p/n)\frac{-1}{z}+\frac{p}{n}
\stesd(z)$. $v_{\esd}$ is the Stieltjes transform of the spectral
distribution of $XX^*/n$.

Currently, the most general version of the result is found in
\citep{silverstein95} and states the following:
\begin{theorem}\label{thm:mapa}
Suppose the data matrix $X$ can be written $X=Y\Sigma_p^{1/2}$,
where $\Sigma_p$ is a $p\times p$ positive definite matrix and $Y$
is an $n\times p$ matrix whose entries are i.i.d (real or complex),
with $E(Y_{i,j})=0$, $E(|Y_{i,j}|^2)=1$ and $E(|Y_{i,j}|^4)<\infty$.

Call $H_p$ the population spectral distribution, i.e the
distribution that puts mass $1/p$ at each of the eigenvalues of the
population covariance matrix, $\Sigma_p$. Assume that $H_p$
converges weakly to a limit denoted $H_{\infty}$. (We write this
convergence $H_p\Rightarrow H_{\infty}$.) Then, when $p,n \tendsto
\infty$, and $p/n\tendsto \gamma$, $\gamma
\in (0,\infty)$,
\begin{enumerate}
\item $v_{\esd}(z)\tendsto \vinfty(z)$, a.s, where $\vinfty(z)$ is a
deterministic  function
\item $\vinfty(z)$ satisfies the equation
\begin{equation}\label{eq:mapa}
-\frac{1}{\vinfty(z)}=z-\gamma\int \frac{\lambda
dH_{\infty}(\lambda)}{1+\lambda \vinfty(z)} \;, \forall z \in
\cplus \tag{M-P}
\end{equation}
\item The previous equation has one and only one solution which is the Stieltjes transform of a measure.
\end{enumerate}
\end{theorem}

In plain English, under the assumptions put forth in Theorem
\ref{thm:mapa}, the spectral distribution of the sample
covariance matrix is asymptotically non-random. Furthermore, it is
fully characterized by the true population spectral distribution,
through the equation (\ref{eq:mapa}).

A particular case of equation (\ref{eq:mapa}) is often of interest:
the situation when all the population eigenvalues are equal to 1.
Then of course, $H_p=\Hinfty=\delta_1$. A little bit of elementary
work leads to the well-known fact in random matrix theory that the
empirical spectral distribution, $\esd$, converges (a.s) to the
\mapa law, whose density is given by, if $\gamma\leq 1$,
$$
f_{\gamma}(x)=\sqrt{(b-x)(x-a)}/(2\pi x\gamma)\;,\;\; \text{ with }
a=(1-\gamma^{1/2})^2\,, b=(1+\gamma^{1/2})^2\;.
$$
We refer the reader to \citep{mp67}, \citep{bai99} and \citep{imj}
for more details and explanations concerning the case $\gamma>1$.
One point of statistical interest is that even though the true
population eigenvalues are all equal to 1, the empirical ones are
now spread on the interval
$[(1-\gamma^{1/2})^2,(1+\gamma^{1/2})^2]$. Plotting the density also
shows that its shape vary with $\gamma$ in a non-trivial way. These
two remarks illustrate some of the difficulties that need to be
overcome when working under ``large $n$, large $p$" asymptotics.

\section{Algorithm and Statistical considerations}\label{sec:Algo}
\subsection{Formulation of the estimation problem}
A remarkable feature of the equation (\ref{eq:mapa}) is that the
knowledge of the limiting distribution of the eigenvalues in the
population given by $\Hinfty$ fully characterizes the limiting
behavior of the eigenvalues of the sample covariance matrix.
However, the relationship between the two is hard to disentangle. As
is common in statistics, the question is how to invert this
relationship to estimate $\Hp$. The question thus becomes, given
$l_1,\ldots,l_p$, the eigenvalues of a sample covariance matrix, can
we estimate the population eigenvalues,
$\lambda_1,\ldots,\lambda_p$, using Equation (\ref{eq:mapa})? Or in
terms of spectral distribution, can we estimate $\Hp$ from $\Fp$?

Our strategy is the following: 1) the first aim is to estimate the
measure $H_{\infty}$ appearing in the \mapa equation. 2) Given an
estimator, $\hatHinfty$, of this measure, we will estimate
$\lambda_i$ as the $i$-th quantile of our estimated distribution. It
is common in statistical practice to get these estimates by using
the $i/(p+1)$ percentile and this is what we do. (We come back to
possible difficulties getting from $\hatHp$ to $\hat{\lambda}_i$ in
\ref{subsubsec:IsoEigen}.) 3) An important point is that since we
are considering fixed distribution asymptotics, our estimate of
$\Hinfty$ will serve as our estimate of $\Hp$, so
$\hatHp=\hatHinfty$.

The main question, then, is how to approach step 1: estimating
$\Hinfty$ based only on $\esd$. Of course, since we can compute the
eigenvalues of $\EmpCovMat$, we can compute $v_{\esd}(z)$ for any
$z$ we choose. By evaluating $v_{\esd}$ at a grid of values
$\{z_j\}_{j=1}^{J_n}$, we have a set of values
$\{v_{\esd}(z_j)\}_{j=1}^{J_n}$ for which equation (\ref{eq:mapa})
should (approximately) hold. We want to find $\hatHinfty$ that will
``best" satisfy equation (\ref{eq:mapa}) across the set of values of
$v_{\esd}(z_j)$. In other words, we will pick
$$
\hatHp=\hatHinfty=\underset{H}{\argmin} \,\,L\left(\left\{\frac{1}{\stesdd(z_j)}+z_j-\frac{p}{n} \int \frac{\lambda
dH(\lambda)}{1+\lambda \stesdd(z_j)}\right\}_{j=1}^{J_n}\right)\;,
$$
where the optimization is over probability measures $H$, and $L$ is
a loss function to be chosen later. In this way we are ``inverting"
the equation (\ref{eq:mapa}), going from $\esd$, an estimate of
$\Finfty$, to an estimate of $\Hinfty$.

We will solve this inverse problem in two steps: discretization and
convex optimization. We give a high-level overview of our method and
postpone implementation details to the Appendix.

To summarize, we face the following interpolation problem: given $J$
an integer and $(z_j,v_{\esd}(z_j))_{j=1}^J$ we want to find an
estimate of $H_{\infty}$ that approximately satisfies equation
(\ref{eq:mapa}). In Section \ref{sec:Consistency}, we show that
doing so for $L_{\infty}$ loss function leads to a consistent
estimator of $\Hinfty$, under the reasonable assumption that all
spectra are bounded.

\subsection{The algorithm}
In order to alleviate the notations, we will replace the notation
$H_{\infty}$ by $H$ when it does not cause any confusion.
\subsubsection{Discretization}\label{subsubsec:discretization}
Naturally, $dH$ can be simply approximated by a weighted sum of
point masses:
$$
dH(x)\simeq \sum_{k=1}^K w_k \delta_{t_k}(x)\;,
$$
where $\{t_k\}_{k=1}^K$ is a grid of points, chosen by us, and
$w_k$'s are weights. The fact that we are looking for a probability
measure imposes the constraints
$$
\sum_{k=1}^K w_k=1\;, \text{ and } w_k\geq 0\;.
$$

This approximation turns the optimization over measures problem into
searching for a vector of weights in $\mathbb{R}_+^K$. After
discretization, the integral in equation (\ref{eq:mapa}) can be
approximated by
$$
\int \frac{\lambda dH(\lambda)}{1+\lambda v}\simeq\sum_{k=1}^K w_k \frac{t_k}{1+t_k v}\;.
$$
Hence finding a measure that approximately satisfies Equation
(\ref{eq:mapa}) is equivalent to finding a set of weights
$\{w_k\}_{k=1}^K$, for which we have
$$
-\frac{1}{\vinfty(z_j)}\simeq z_j-\frac{p}{n}\sum_{k=1}^K w_k
\frac{t_k}{1+t_k\vinfty(z_j)} \;, \forall j\;.
$$

Naturally, we do not get to observe $\vinfty$, and so we make a
further approximation and replace $\vinfty$ by $v_{\esd}$. Our
problem is thus to find $\{w_k\}_{k=1}^K$ such that
$$
-\frac{1}{v_{\esd}(z_j)}\simeq z_j-\frac{p}{n}\sum_{k=1}^K w_k
\frac{t_k}{1+t_k v_{\esd}(z_j)} \;, \forall j\;.
$$

 One good thing
about this approach is that the problem we now face is linear in the
weights, which are the only unknowns here. We will demonstrate that
this allows us to cast the problem as a relatively simple convex
optimization problem.
\subsubsection{Convex Optimization
formulation}\label{subsubsec:convexOptFormulation} To show that we
can formulate our inverse problem as a convex problem, let us call
the approximation errors we make
$$
e_j=\frac{1}{v_{\esd}(z_j)}+z_j-\frac{p}{n}\sum_{k=1}^K w_k
\frac{t_k}{1+v_{\esd}(z_j) t_k}\;.
$$
As explained above, there are two sources of error in $e_j$: one
comes from the discretization of the integral involving $\Hinfty$.
The other one comes from the substitution of $\vinfty$, a non-random
and asymptotic quantity, by $v_{\esd}$, a (random) quantity
computable from the data. $e_j$ is of course a complex number in
general.

We can now state several convex problems as approximation of the
inversion of the \mapa equation problem. We show in Section
\ref{sec:Consistency} consistency of the solution of the ``$L_{\infty}"$ version of the
problem described below. Here are a few examples of convex
formulations for our inverse problem. In all these problems, the
$w_k$'s are constrained to sum to 1 and to be non-negative.
\begin{enumerate}
\item ``$L_{\infty}$" version: Find $w_k$'s to
$$
\text{Minimize } \max_{j=1,\ldots,J_n} \max
\left\{\left|\myreal{e_j}\right|,\left|\imag{e_j}\right|\right\}
$$
\item ``$L_{2}$" version: Find $w_k$'s to
$$
\text{Minimize } \sum_{j=1}^{J_n} \left|e_j\right|\;.
$$
\item ``$L_{2}$-squared" version: Find $w_k$'s to
$$
\text{Minimize } \sum_{j=1}^{J_n} \left|e_j\right|^2\;.
$$
\end{enumerate}
The advantages of formulating our problem as a convex optimization
problem are many. We will come back to the more statistical issues
later. From a purely numerical point of view, we are guaranteed that
an optimum exists, and fast algorithms are available. In practice,
we used the optimization package MOSEK (see \citep{mosek}), within
Matlab, for solving our problems.

Because the rest of the article focuses particularly on the
``$L_{\infty}$" version of the problem described above, we want to
give a bit more details about it. The ``translation" of the problem
into a convex optimization problem is
\begin{gather*}
\min_{(w_1,\ldots,w_K,u)} u\\
\forall j,\; -u\leq \myreal{e_j} \leq u\\
\forall j,\; -u\leq \imag{e_j} \leq u\\
\text{ subject to }
\sum_{i=1}^K w_k =1\\
\text{ and } w_k \geq 0, \forall k
\end{gather*}
This is a linear program (LP) with unknowns $(w_1,\ldots,w_K)$ and
$u$ (see \citep{boydvandenberghe04} for standard manipulations to
make it a standard form LP).

The simulations we present in Section \ref{sec:Simulations} were
made using this version of this algorithm. The proof in Section
\ref{sec:Consistency} applies to this version of the algorithm.

\subsection{Statistical considerations}
The formulation we proposed is quite flexible and has several
important qualities. For instance, regularization constraints can be
easily handled through our proposal. We also can view the algorithm
as a form of ``basis pursuit" in measure space, from which we can
draw some practical conclusions.

\subsubsection{Regularization and constraints}
Methods to invert the \mapa equation should be flexible enough to
accommodate reasonable constraints that could provide additional
improvement to our estimate of $\Hp$. The fact that we essentially
just optimize over the weights $w_k$'s mean that we can easily
regularize and add constraints. For instance, we might want to
regularize our estimator and make it smoother by adding a total
variation penalty (on the $w_k$'s) to our objective function. In
terms of constraints, we might want to specify that the first moment
of our estimate $\hatHp$ match the trace of $\EmpCovMat/p$, since we
know that the trace of $\EmpCovMat/p$ is a good estimate of the
trace of $\Sigma_p/p$ (see e.g
\citep{jonsson82}), and that the trace of $\Sigma_p/p$ is equal to the
first moment of $\Hp$. Note that constraints on the moments of our
estimator are linear in the $w_k$'s and so such constraints would
still lead to a convex problem. The framework we provide can very
easily incorporate these two examples of penalty and constraints, as
well as many others.

\subsubsection{A ``basis pursuit" point of view}
A semantic point is needed before we start our discussion. We use
the term ``basis pursuit" in a loose sense: we are not referring to
the algorithm proposed in
\citep{chendonohosaunders98} but rather use this expression as a generic term for describing
techniques that aim to optimize the representations of functional
objects in overcomplete dictionaries. We refer the reader to
\citep[Chapter 5]{htf01} for some of the core statistical ideas of
these so-called basis expansion methods.

The algorithm we propose can be viewed as a relaxation of a measure
estimation problem. We want to estimate a measure $\Hinfty$ and
instead of searching among all possible probability measures, we
restrict our search space to mixtures of certain class of
probability measures. In
\ref{subsubsec:discretization} for instance, we restricted the
choice to mixture of point masses. In that sense, we can view it as
a type of ``basis pursuit" in probability measure space. We first
choose a ``dictionary" of probability measures on the real line, and
we then decompose our estimator on this dictionary, searching for
the best coefficients. Hence our problem can be formulated as
$$
\text{ find the best possible weights } \{w_1,\ldots,w_N\} \text{ with }d\widehat{H}
=\sum_{i=1}^N w_i dM_i
$$
where the $M_i$'s are the measures in our dictionary.

In the preceding discussion on discretization, we restricted
ourselves to $M_i$'s being point masses at chosen ``grid points". Of
course, we can enlarge our dictionary to include, for instance:
\begin{enumerate}
\item  Probability measures that are uniform on an interval:
$dM_i(x)=1_{x\in [a_i,b_i]}dx/(b_i-a_i)$.
\item  Probability measures that have a linearly increasing density on an interval $[a_i,b_i]$
and density $0$ elsewhere. So
$dM_i(x)=1_{[a_i,b_i]}2(x-a_i)/(b_i-a_i)^2 dx$, and density $0$
elsewhere.
\item  Probability measures that have a linearly decreasing density on an interval $[a_i,b_i]$,
and density $0$ elsewhere. So
$dM_i(x)=1_{[a_i,b_i]}2(b_i-x)/(b_i-a_i)^2 dx$.
\end{enumerate}

If we decide to include a probability measure $M$ in our dictionary,
the only requirement is that we be able to compute the integral
$$
\int \frac{\lambda dM(\lambda)}{1+\lambda v}
$$
for any $v$ in $\cplus$.

Choosing a larger dictionary increases the size of the convex
optimization problems we try to solve, and hence is at first glance
computationally harder. However, statistically, enlarging the
dictionary may lead to sparser representations of the measure we are
estimating, and hence, at least intuitively, lead to better
estimates of $\Hinfty$. The most favorable case is of course when
$\Hinfty$ is a mixture of a small number of measures present in our
dictionary. For instance, if $\Hinfty$ has a density whose graph is
a triangle, having measures as described in points 2 and 3 above
would most likely lead to sparser and maybe more accurate estimates.
In the presence of a priori information on $\Hinfty$, the choice of
dictionary should be adapted so that $\Hinfty$ has a sparse
representation in the dictionary.
\subsubsection{Useful properties of the algorithm}
One important advantage of choosing to estimate measures instead of
choosing to estimate a high-dimensional vector is that the
algorithm's complexity does not increase with the size of the answer
required by the user. Hence given a $p$ dimensional vector of
eigenvalues, once the values $v_{\esd}(z_j)$ are computed, the
computational cost of the algorithm is the same irrespective of $p$.
This means that for large $p$ problems, only one difficult
computation is required: that of the eigenvalues of the empirical
covariance matrix. Our algorithm is hence, in some sense,
``dimension-free", i.e, except for the computation of the
eigenvalues, it is insensitive to the dimensionality of our original
problem. This scaling property is important for high-dimensional
problems.

Another good property of our method is that it is independent of the
basis in which the data is represented. Because our method requires
only as input the eigenvalues of the sample covariance matrix -
quantities obviously independent of the original basis of the data -
our method is basis independent.

In other respects, Theorem \ref{thm:mapa} holds for random variables
that have a 4-th moment; we are not limited to Gaussian random
variables. Complex random variables are also possible. Hence, the
theorem is well-suited for wide applicability. Elementary properties
of Gaussian random variables show that Theorem
\ref{thm:mapa} covers all possible Gaussian problems. This will not
be true for all distributions, but the scope of the theorem is still
very wide. Note also that the Equation (\ref{eq:mapa}) holds in
greater generality than mentioned in Theorem
\ref{thm:mapa}. We refer the reader to the original paper
\citep{mp67} for further examples, in particular when the data is
distributed on spheres or ellipsoids. (The original formulation of
the theorem allows for dependence between the entries of the matrix
$Y$, but the convergence is not shown to be almost sure.)

\subsubsection{The case $p>n$ and how large is large?}
Another advantage of the proposed method is that it is insensitive
to whether $p$ is larger than $n$ or $n$ is larger than $p$. The
only requirement is that they both be quite large. We had reasonable
to good results in simulation as soon as $p>30$ or so. As a matter
of fact, it is quite clear that to have reasonably accurate
estimates of the eigenvalues, we need to ``populate" the interval
$[\lambda_p,\lambda_1]$ with enough points, for otherwise quantile
methods may be somewhat inaccurate.

\subsubsection{On covariance estimation, linear and non-linear shrinkage of eigenvalues}

There is some classical and more recent statistical work on
shrinkage of eigenvalues to improve covariance estimation. We refer
the reader to Section 4.1 in \citep{LedoitWolf04} for some examples
due to Charles Stein and Leonard Haff, unfortunately in unpublished
manuscripts. More recently, in the interesting paper by
\citep{LedoitWolf04}, what was proposed is to linearly shrink the
eigenvalues of $\EmpCovMat$ toward the identity : i.e $l_i$'s become
$\tilde{l}_i=(1-\rho)l_i+\rho$'s, for some $\rho$, independent of
$i$, chosen using the data and the
\mapa law. Then the authors of
\citep{LedoitWolf04} proposed to estimate $\Sigma_p$ by $(1-\rho)\EmpCovMat+\rho
Id_p$. Since this latter matrix and $\EmpCovMat$ have the same
eigenvectors, their method of covariance estimation can be viewed as
linearly shrinking the sample eigenvalues and keeping the
eigenvectors of $\EmpCovMat$ as estimates of the eigenvectors of
$\Sigma_p$.

Our method of estimation of the population eigenvalues can be viewed
as doing a non-linear shrinkage of the sample eigenvalues. While we
could propose to just keep the eigenvectors of $\EmpCovMat$ as
estimates of the eigenvectors of $\Sigma_p$, and hence get an
estimate of the population covariance matrix, we think one should be
able to do better by using the eigenvalue information to drive the
eigenvector estimation. It is known that in ``large $n$, large $p$"
asymptotics, the eigenvectors of the sample covariance matrix are
not consistent estimators of the population eigenvectors (see
\citep{debashis}), even in the most favorable cases. However,
having a good idea of the structure of the population eigenvalues
should help us estimate the eigenvectors of the population
covariance matrix, or at least formulate the right questions for the
problem at hand. For instance, the inferred structure of the
covariance matrix could help us decide how many subspaces we need to
identify: if, for example,  it turned out that the population
eigenvalues were clustered around two values, we would have to
identify two subspaces, the dimensions of these subspaces being the
number of eigenvalues clustered around each value. Also, having
estimates of the eigenvalues tell us how much variance our
``eigenvectors" will have to explain. In other words, our hope is
that taking advantage of the crucial eigenvalue information we are
now able to gather will lead to better estimation of $\Sigma_p$ by
doing a ``reasoned" spectral decomposition. Work in this direction
is in progress.

\subsubsection{Asymptotics at fixed spectral distribution and isolated eigenvalues}\label{subsubsec:IsoEigen}

Our algorithm actually uses asymptotics assuming a fixed spectral
distribution: we are essentially fixing $H_p=\Hinfty$ when solving
our optimization problem. Naturally, this does not mean that $p$ is
fixed. Note that this is what is classically done is statistics: for
the simple problem of estimating the mean of a population from a
sample $Z_1,\ldots,Z_K,$ it is common to assume that the $Z_k$'s
have the same mean $\mu$, and that $\mu$ does not depend on $K$.
However, when studying the asymptotic properties of this simple
estimator, we could require to actually have $\mu(K)$, with
$\mu(K)\tendsto
\mu$. (All we would have to do is have a triangular array of data,
and getting to observe just one row of this array at a time.) Hence
our fixed spectral distribution ``assumption" is very natural and
similar to classical assumptions made in estimation problems.

Let us go back now to the problem of isolated eigenvalues. Suppose
we get to see data in $\mathbb{R}^{p_0}$ for some $p_0$. Then, any
isolated eigenvalue that may be present is numerically treated as if
the mass that is attached to it is held fixed at $1/p_0$ when
$p\tendsto
\infty$. So a point mass at the
corresponding population eigenvalue would appear in $\hatHp$. This
has been verified numerically. If the estimator were perfect, this
mass should be equal to $1/p_0$. However, because of variability it
may not be exactly of mass $1/p_0$. Then, estimating the population
eigenvalues by the quantiles of the estimated population spectral
distribution, we may ``miss" this isolated eigenvalue. In the case
of the largest eigenvalue, that would happen if the mass found
numerically at this isolated eigenvalue is less than $1/(p_0+1)$. So
isolated eigenvalues will require special care and caution,
particularly in going from $\hatHp$ to $\hat{\lambda_i}$. While the
method focuses on identifying the structure of the population
eigenvalues and hence may have problems when it comes to estimating
isolated eigenvalues, we have found in practice that it still
provided a good tool for this task but that some care was required.

\subsubsection{Existing related work}
\label{subsubsec:ChooseEstimMeas}
As far as we know, there has been no work on non-parametric
estimation of $\Hp$ or $\Hinfty$ using the \mapa equation. However,
some work exists in the Physics' literature
(\citep{burdaetal04,burdajurkwaclaw}), that takes advantage of the
\mapa law to estimate some moments of $\Hinfty$. $\Hinfty$ is then
assumed to a be a mixture of a finite and pre-specified number of
point masses (see \citep[p. 303]{burdaetal04}) and the moments are
then matched with possible point masses and weights. While these
methods might be of some use sometimes, we think they require too
many assumptions to be practically acceptable for a broad class of
problems. It might be tempting to try to develop an non-parametric
estimator from moments, but we think that without the strong
assumptions made in \citep{burdaetal04}, those estimators will
suffer drastically from: 1) the number of moments needed a priori
may be large, and large moments are very unreliable estimators; 2)
moments estimated indirectly may not constitute a genuine family of
moments: certain Hankel matrices need to be positive semi-definite
and will not necessarily be so. Semi-definite programming type
corrections will then be necessary, but hard to implement. 3) Even
if one has a genuine moment sequence, there are usually many
distributions with the same moments. Choosing between them is
clearly going to be a difficult task.

\section{Simulations}\label{sec:Simulations}
We now present some simulations to illustrate the practical
capabilities of the method. The objectives of eigenvalues estimation
are many-folds and depend of the area of applications. We review
some of those that inspired our work.

In settings like PCA, one basically wishes to discover some form of
structure in the covariance matrix by looking at the eigenvalues of
the sample covariance matrix. In particular, a situation where the
population eigenvalues are different from each other indicates that
projecting the data in some projections will be more ``informative"
that projecting it in other directions; while in the case where all
the population eigenvalues are equal, all projections are equally
informative or uninformative. As our brief discussion of the
\mapa law illustrated, in the ``large $n$, large
$p$" setting, it is difficult to know from the sample eigenvalues
whether all population eigenvalues are equal to each other or not,
or even if there is any kind of structure in them. When $p$ and $n$
are both large, standard graphical methods like the scree plot tend
to look similar whether or not there is structure in the data. We
will see that our approach is able to differentiate between the
situations. Among other things, our method can thus be thought as a
alternative to the scree plot for high-dimensional problems.

In other applications, one focuses more on trying to estimate the
value of the largest or smallest eigenvalues. In PCA, the largest
population eigenvalues measure how much variance we can explain
through a low dimensional projection and is hence important. In
financial applications, like the Markovitz' portfolio optimization
problem, the small population eigenvalues are important. They
essentially measure what is the minimum risk one can take by
investing in a portfolio of certain stocks (see \citep{lalouxetal}
and
\citep[Chapter 5]{CampbellLoMacKinlay}). However, as explained in the Appendix, the largest eigenvalue of
the sample covariance matrix tends to overestimate the largest
eigenvalue of the population covariance. And similarly, the smallest
eigenvalue of the sample covariance matrix tends to underestimate
its population counterpart. What that means is that using these
measures of ``information" and ``risk", we will tend to overestimate
the amount of information there is in our data and tend to
underestimate the amount of risk there is in our portfolios. So it
is important to have tools to correct this bias. Our estimator
provides a way to do so.

\subsection{Details of the simulations}
We illustrate the performance of our method on three cases, each
with very different covariance structure. We will give more details
on each individual case in the following subsections.

We now describe more precisely these examples. The first case is
that of $\Sigma_p=\id_p$, in other words, there is no ``information"
in the data. However standard graphical statistical methods like the
``scree plot" will tend to show a pattern in the eigenvalues. We
will show that our method is generally able to inform us that all
the eigenvalues are equal.

The second  case is one where $\Sigma_p$ has 50\% of its eigenvalues
equal to 1 and 50\% equal to 2. While it should be easy to discern
that there are two very distinct clusters of eigenvalues in the
population, in high-dimension the sample eigenvalues will often blur
the clusters together. We show that our method generally recovers
these two clusters well.

Finally, the third example is one where $\Sigma_p$ is a Toeplitz
matrix. More details on Toeplitz matrices are given in
\ref{subsubsec:ToepCov}. This situation poses a
harder estimation problem. While the asymptotic behavior of the
eigenvalues of such matrices is well understood, there are generally
no easy and explicit formulas to represent the limit. We present the
results to show that even in this difficult setting, our method
performs quite well.

To measure the performance of our estimators, we compare the L\'evy
distances between our estimator, $\hatHp$, and the true distribution
of the population eigenvalues, $\Hp$, to that of the empirical
spectral distribution, $\Fp$, to $\Hp$. Our choice is motivated by
the fact that the L\'evy distance can be used as a metric for weak
convergence of distributions on $\mathbb{R}$. Recall (see e.g
\citep{durrett96}) that the L\'evy distance between two distributions $F$ and $G$
on the real line is defined as
$$
d_L(F,G)=\inf\{\eps>0: F(x-\eps)-\eps \leq G(x)\leq F(x+\eps)+\eps
\,
\;,\forall x\}\;.
$$

In the plots we will depict the cumulative distribution function
(cdf) of our estimated measures. Recall that the estimates of the
population eigenvalues $\lambda_i$'s are obtained by taking
appropriate percentiles of these measures.
\subsubsection{The case $\Sigma_p=\id_p$}
\begin{figure}
\centering
\subfloat[Eigenvalues (scree plot) of the sample covariance matrix]{\label{fig:Id500100Case-screeplot}\includegraphics[scale=.3,angle=-90]{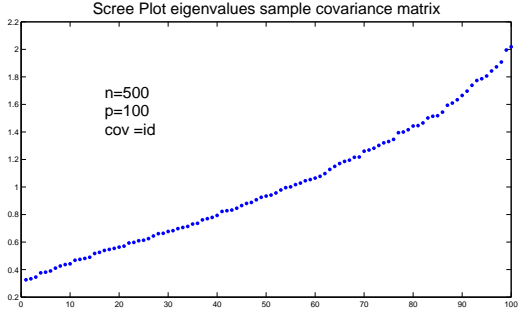}}\qquad
\subfloat[CDF eigenvalues, sample covariance matrix ($\Fp$)]{\includegraphics[scale=.3,angle=-90]{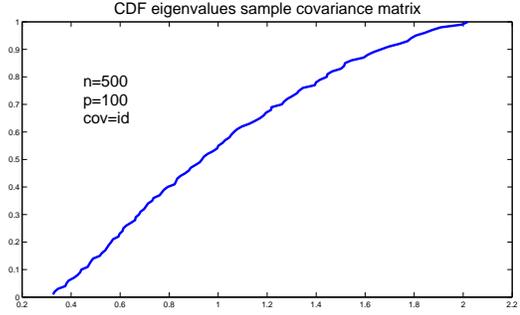}}\qquad
\subfloat[CDF eigenvalues, estimated population covariance matrix ($\hatHp$)]{\label{fig:Id500100Case-mymethod}\includegraphics[scale=.3,angle=-90]{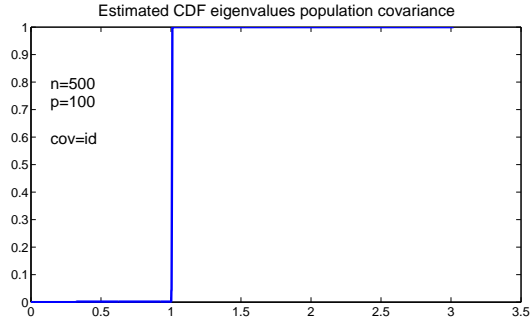}}\qquad
\caption[Illustration of the performance of the estimator:]{case $\boldsymbol{\Sigma_p=\id_p}$.
The three figures above compare the performance of our estimator to
the one derived from the sample covariance matrix on one realization
of the data. The data matrix $X$ is $500\times 100$. All its entries
are iid ${\cal N}(0,1)$. The population covariance is
$\Sigma_p=\id_{100}$, so the distribution of the eigenvalues is a
point mass at 1. This is what our estimator (Figure
\subref{fig:Id500100Case-mymethod}) recovers. Average computation
time (over 1000 repetitions) was 13.33 seconds, according to
\texttt{Matlab} tic and toc functions. Implementation
details
are in the Appendix.}
\label{fig:Id500100Case}
\end{figure}

In this situation, the \mapa law predicts that instead of being
concentrated at 1 like the population eigenvalues, the sample
eigenvalues will be spread on the interval
$[(1-\sqrt{p/n})^2,(1+\sqrt{p/n})^2]$. This is problematic, since by
looking at the scree plot of just the sample eigenvalues, one might
think that some population eigenvalues are (much) larger than others
and hence some projections of the data are more informative than
others. This is vividly illustrated on Figure
\ref{fig:Id500100Case-screeplot}. However, as we see on Figure
\ref{fig:Id500100Case-mymethod}, the method we propose finds that
the population spectral distribution is very close to a point mass
at 1, and all eigenvalues are thus close to 1. Statistically, this
of course means that there is no preferred direction to project the
data. All directions are equally informative, or uninformative.

\begin{figure}
\begin{center}
\includegraphics[scale=.25]{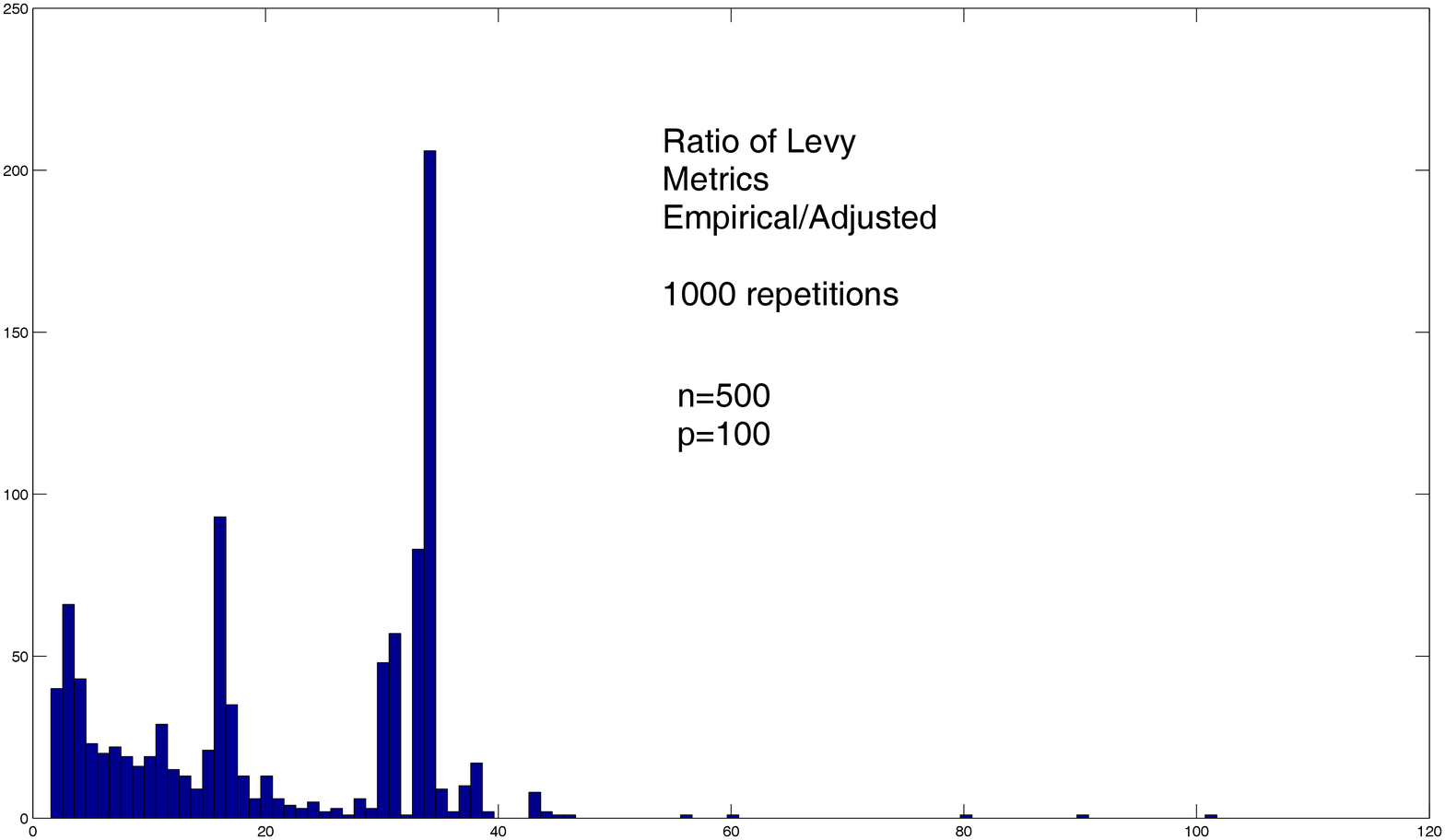}
\end{center}
\caption{\label{fig:ratioLevyDIdCase}case
$\boldsymbol{\Sigma_p=\id_p}$: Ratios $d_L(\hatHp,\Hp)/d_L(\Fp,\Hp)$
over 1,000 repetitions. Dictionary consisted of only point masses.
Large values indicate better performance of our algorithm. All
ratios were found to be larger than 1. }
\end{figure}

The figures presented in Figure \ref{fig:Id500100Case} were chosen
at random among 1000 Monte-Carlo simulations and are very
encouraging. To further our empirical investigation of the
performance of our method, we repeated the estimation process 1000
times. Another advantage is that on further investigation (manually
checking the graphs of many of the estimators we obtained) we saw
that the estimator consistently gets the structure ``right", namely
a huge spike in the vicinity of 1. This is of course very important
for applications such as PCA, where the structure of the spectrum of
the covariance matrix is of fundamental importance. For each
repetition, we estimated the distribution of the eigenvalues in the
population, and computed the L\'evy distance of our estimator,
$\hatHp$, to the true distribution, $\Hp$, in this case a point mass
at 1. We did the same for the empirical spectral distribution $\Fp$.
Figure
\ref{fig:ratioLevyDIdCase} shows the ratio
$d_L(\hatHp,\Hp)/d_L(\Fp,\Hp)$ for these simulations. Our estimator
clearly outperforms the one derived from the sample covariance
matrix, often by a dramatic factor.

\subsubsection{The case $H_p=.5 \delta_1 + .5 \delta_2$}
In this case the eigenvalues of the population covariance matrix are
split into two clusters of equal size. For the specific example we
investigate, 50\% of the eigenvalues are equal to 1 and 50\% are
equal to 2.

While it should be easy to discern that there are two very distinct
clusters of population eigenvalues, when $p$ is sufficiently close
to $n$ the two clusters merge together and the scree plot of the
sample eigenvalues does not show a clear separation between the two
regions. The
\mapa law predicts (in the case of identity covariance) that the
sample eigenvalues spread over larger and larger intervals as $p$
gets closer to $n$. Therefore, it is intuitively not surprising that
when we have two not too distant clusters of population eigenvalues,
the corresponding sample eigenvalues  would start to overlap if $p$
is close enough to $n$.

\begin{figure}
\centering
\subfloat[Scree plot of eigenvalues, sample covariance matrix: no clear separation around the 50th eigenvalue]{\label{fig:TwoPoints500100Case-screeplot}\includegraphics[scale=.3,angle=-90]{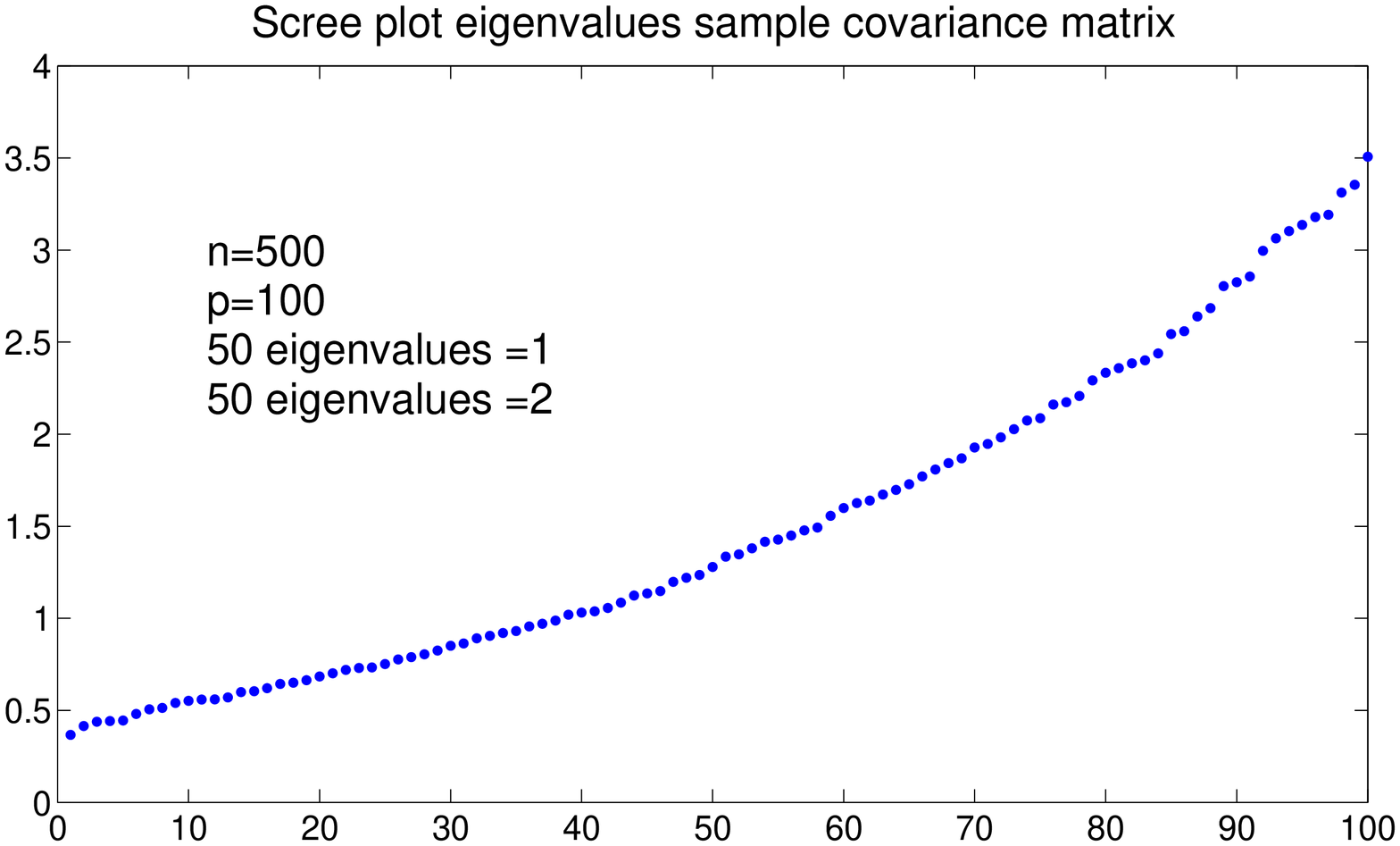}}\qquad
\subfloat[CDF eigenvalues sample covariance matrix ($\Fp$)]{\includegraphics[scale=.3,angle=-90]{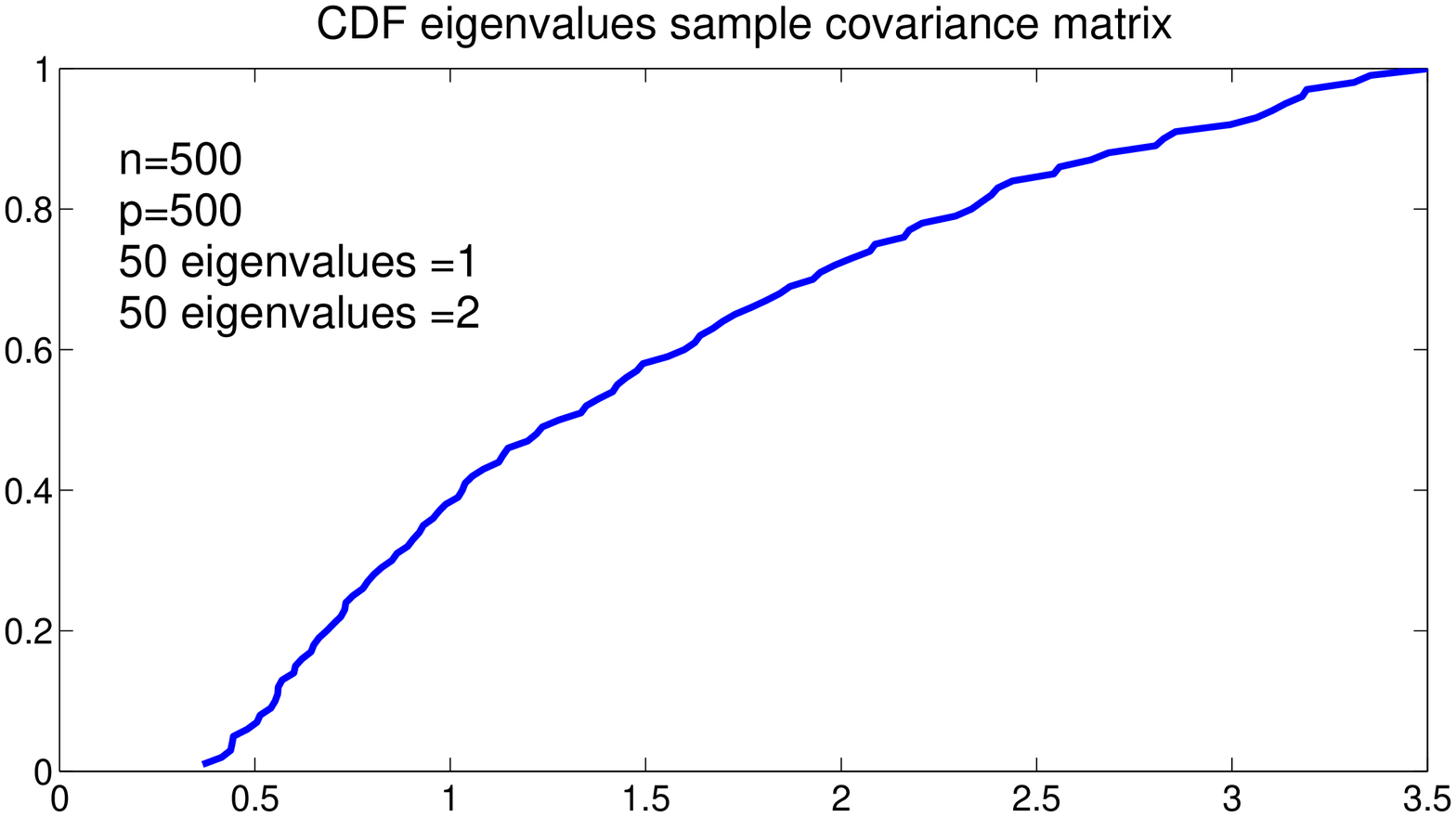}}\qquad
\subfloat[Estimated CDF of eigenvalues of
population covariance matrix
($\hatHp$)]{\label{fig:TwoPoints500100Case-mymethod}\includegraphics[scale=.3,angle=-90]{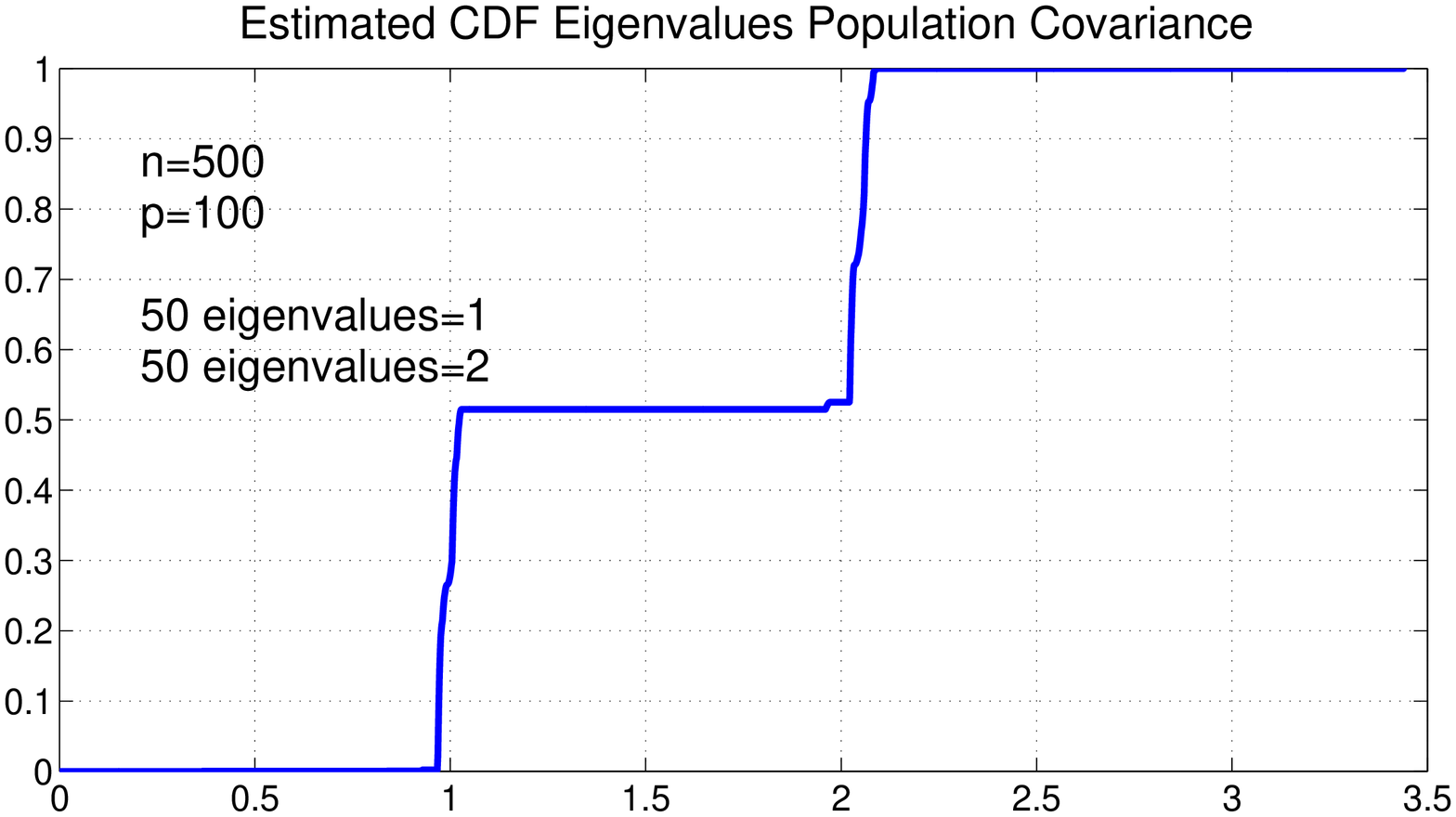}}\qquad
\caption[Illustration of the estimator]{case
$\boldsymbol{H_p=.5 \delta_1 + .5 \delta_2}$ : the three figures
above compare the performance of our estimator on one realization of
the data. The data matrix $Y$ is $500\times 100$. All its entries
are iid ${\cal N}(0,1)$. The covariance is diagonal and has spectral
distribution $H_p=.5 \delta_1 + .5
\delta_2$. In other words, 50 eigenvalues are equal to 1 and fifty
eigenvalues are equal to 2. This is essentially what our estimator
(Figure
\subref{fig:TwoPoints500100Case-mymethod}) recovers. Average
computation time (over 1000 repetitions) was 15.71 seconds,
according to \texttt{Matlab} tic and toc functions.} 
\label{fig:TwoPoints500100Case}
\end{figure}

We did a Monte Carlo analysis (similar to the one done in the case
of $\id_p$ covariance) of our estimator and did comparisons to the
empirical spectral distribution. As in the case of $\id_p$, we
present a figure showing the ratio of the L\'evy distance of the two
estimates to the true distribution.
\begin{figure}
\begin{center}
\includegraphics[scale=.25]{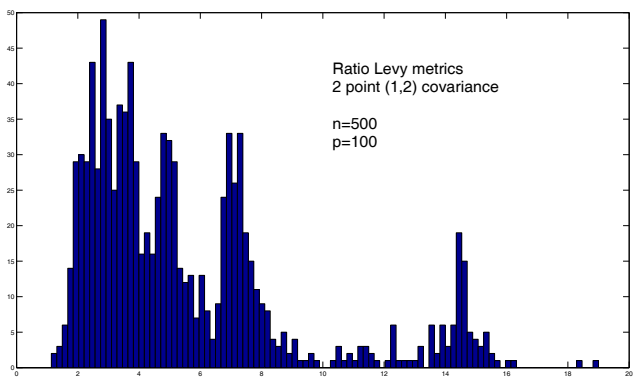}
\end{center}
\caption{\label{fig:ratioLevyD2Point}case
$\boldsymbol{H_p=.5\delta_1 + .5 \delta_2}$ : Ratios
$d_L(\hatHp,\Hp)/d_L(\Fp,\Hp)$ over 1,000 repetitions. Dictionary
consisted of only point masses. Large values indicate better
performance of our algorithm. All ratios were found to be larger
than 1.}
\end{figure}
Figure \ref{fig:ratioLevyD2Point} shows that once again our
estimator clearly outperforms the one derived from the sample
covariance matrix, by a large factor. Again, upon further
investigation, the estimator generally gets the correct structure of
the distribution of the population eigenvalues: in this case two
spikes at 1 and 2.

\subsubsection{The case of a Toeplitz covariance
matrix}\label{subsubsec:ToepCov} Finally, we performed the same type
of analysis on a Toeplitz matrix, to show that the method we propose
works quite well on more complicated types of covariance structures.
Note that generally this is inherently a quite difficult problem, if
we do not assume a priori that we know that the matrix is Toeplitz.

We recall that a Toeplitz matrix $T$ is a matrix whose entries
satisfy $T_{i,j}=t(i-j)$, for a certain function $t$. Since
covariance matrices are symmetric, the Toeplitz matrices at hand
will satisfy $T_{i,j}=t(|i-j|)$. The limiting spectral distribution
of these objects are very well understood: see
\citep{bottchersilbermann}, \citep{gray} or \citep{grenanderszego58}.

Approaches exist that take advantage of the particular structure of
a Toeplitz matrix. See for instance, the interesting papers
\citep{bickellevina04} and for even more generality - beyond Toeplitz
matrices -
\citep{bickellevina06}. However, these approaches are very basis
dependent; they assume that the variables are measured in the
appropriate basis.  In data analysis, this may sometimes be
justified and sometimes not. In particular, if the order of the
variables is permuted, the resulting estimators might change. Since
we want to be able to avoid this type of behavior, we feel that a
``basis independent" method is needed and should be available.
Finding such a method was one of the original motivations of our
investigations.

\begin{figure}
\centering
\subfloat[Scree plot, Eigenvalues sample covariance matrix]{\label{fig:Toep500100Case-screeplot}\includegraphics[scale=.3,angle=-90]{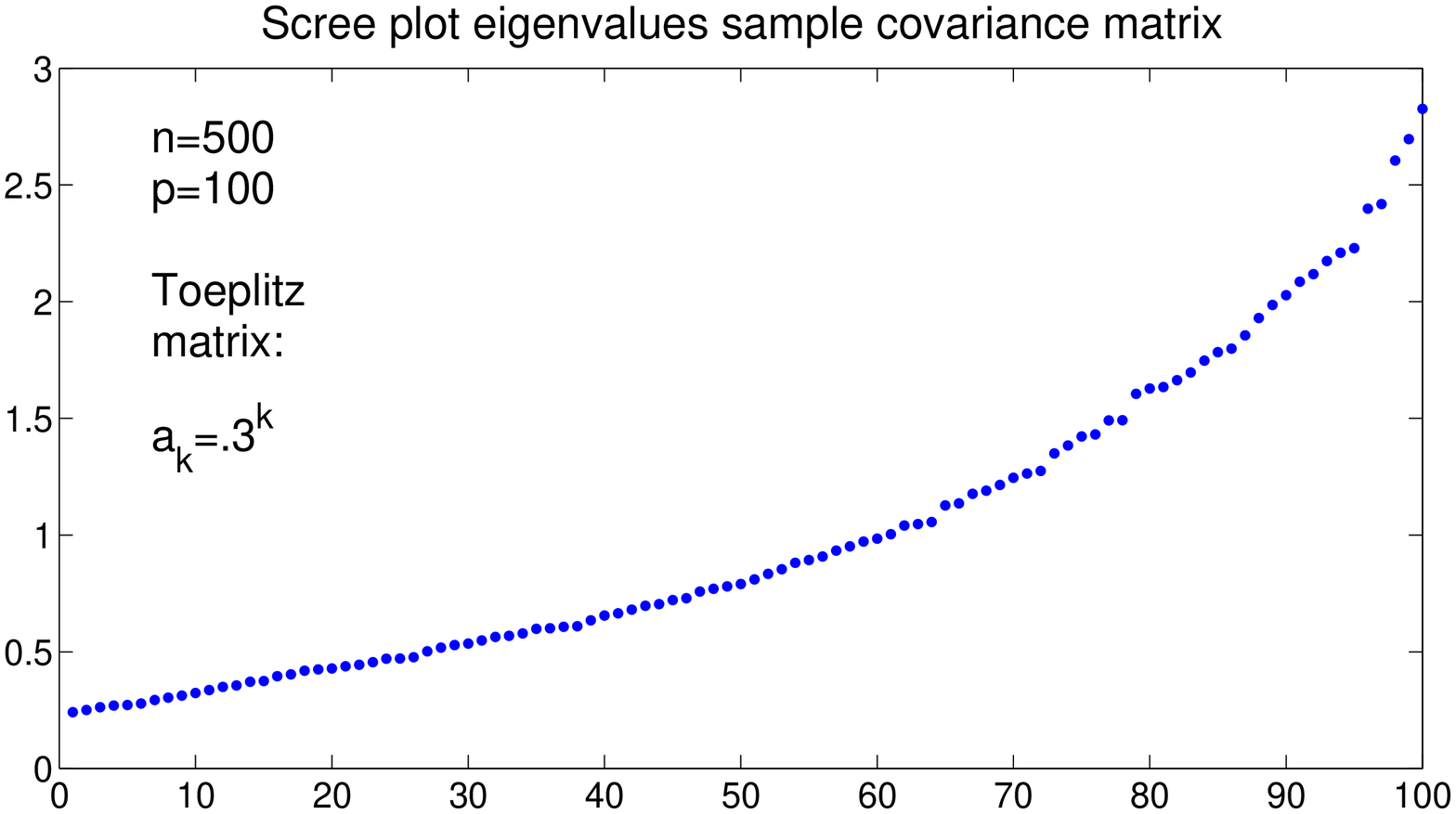}}\qquad
\subfloat[CDF eigenvalues sample covariance matrix ($\Fp$)]{\includegraphics[scale=.3,angle=-90]{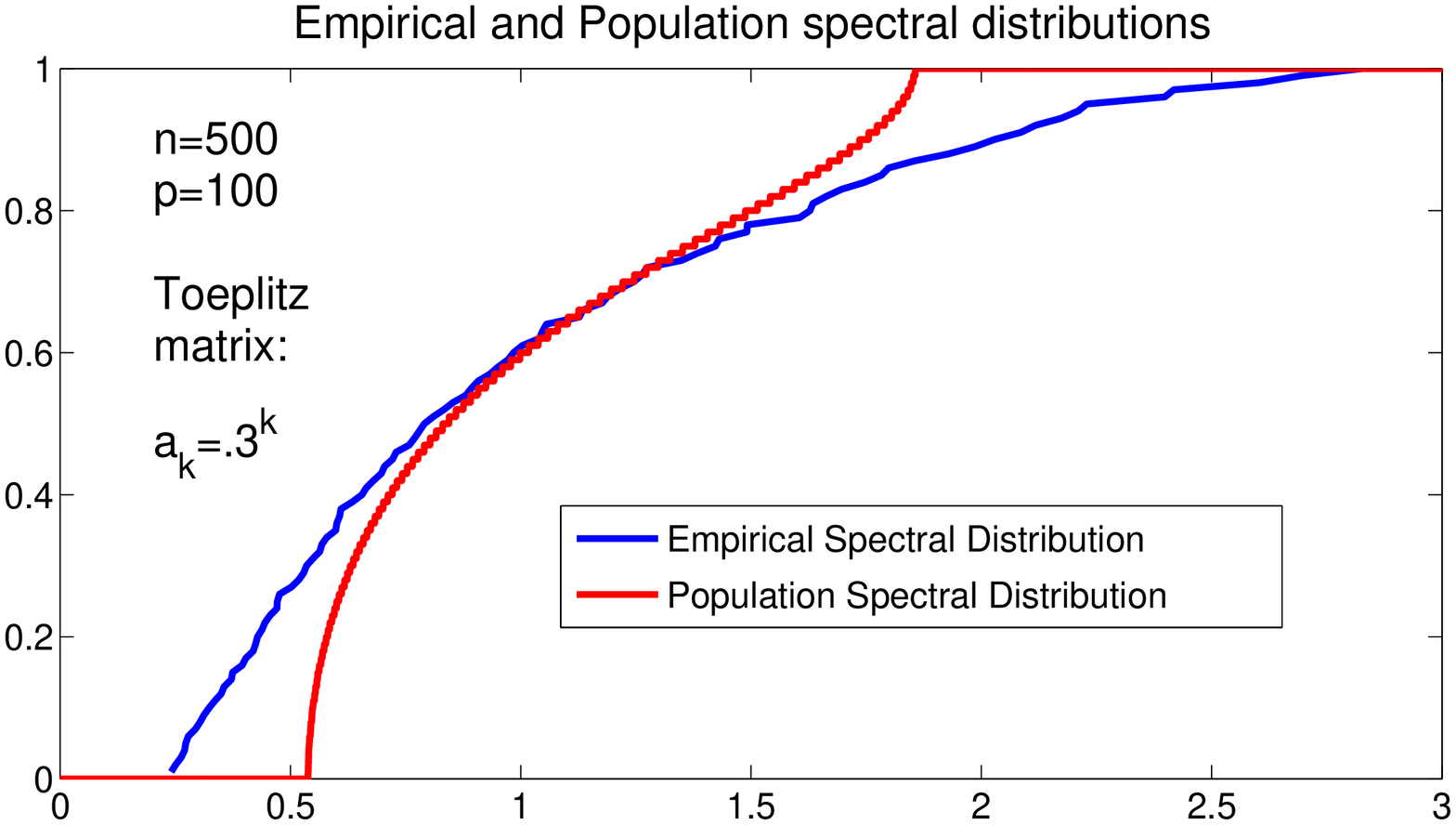}}\qquad
\subfloat[Estimated CDF of eigenvalues of population
covariance matrix
($\hatHp$)]{\label{fig:Toep500100Case-mymethod}\includegraphics[scale=.3,angle=-90]{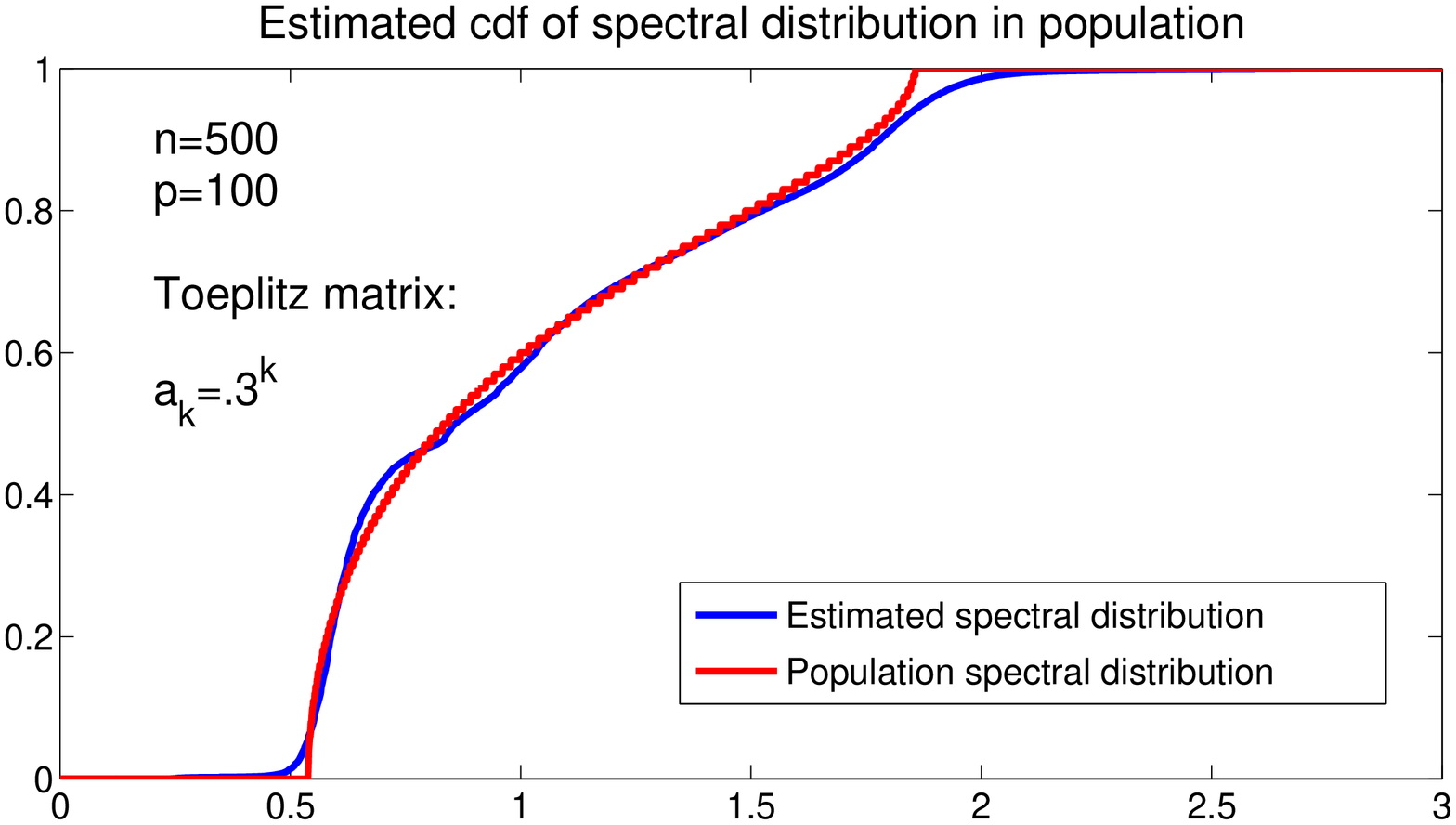}}\qquad
\caption[Illustration of the estimator]{case
$\boldsymbol{\Sigma_p}$ \textbf{Toeplitz with entries}$\boldsymbol{
.3^{|i-j|}}$ : the three figures above show the performance of our
estimator on one realization of the data. The data matrix $Y$ is
$500\times 100$. All its entries are iid ${\cal N}(0,1)$. The
covariance is Toeplitz, with $t(|i-j|)=.3^{|i-j|}$. In Figure
\subref{fig:TwoPoints500100Case-mymethod}, we superimpose
our estimator (blue curve) and the true distribution of eigenvalues
(red curve). Average computation time (over 1000 repetitions) was
16.61 seconds, according to \texttt{Matlab} tic and toc functions.}
\label{fig:Toep500100Case}
\end{figure}

\begin{figure}
\begin{center}
\includegraphics[scale=.25]{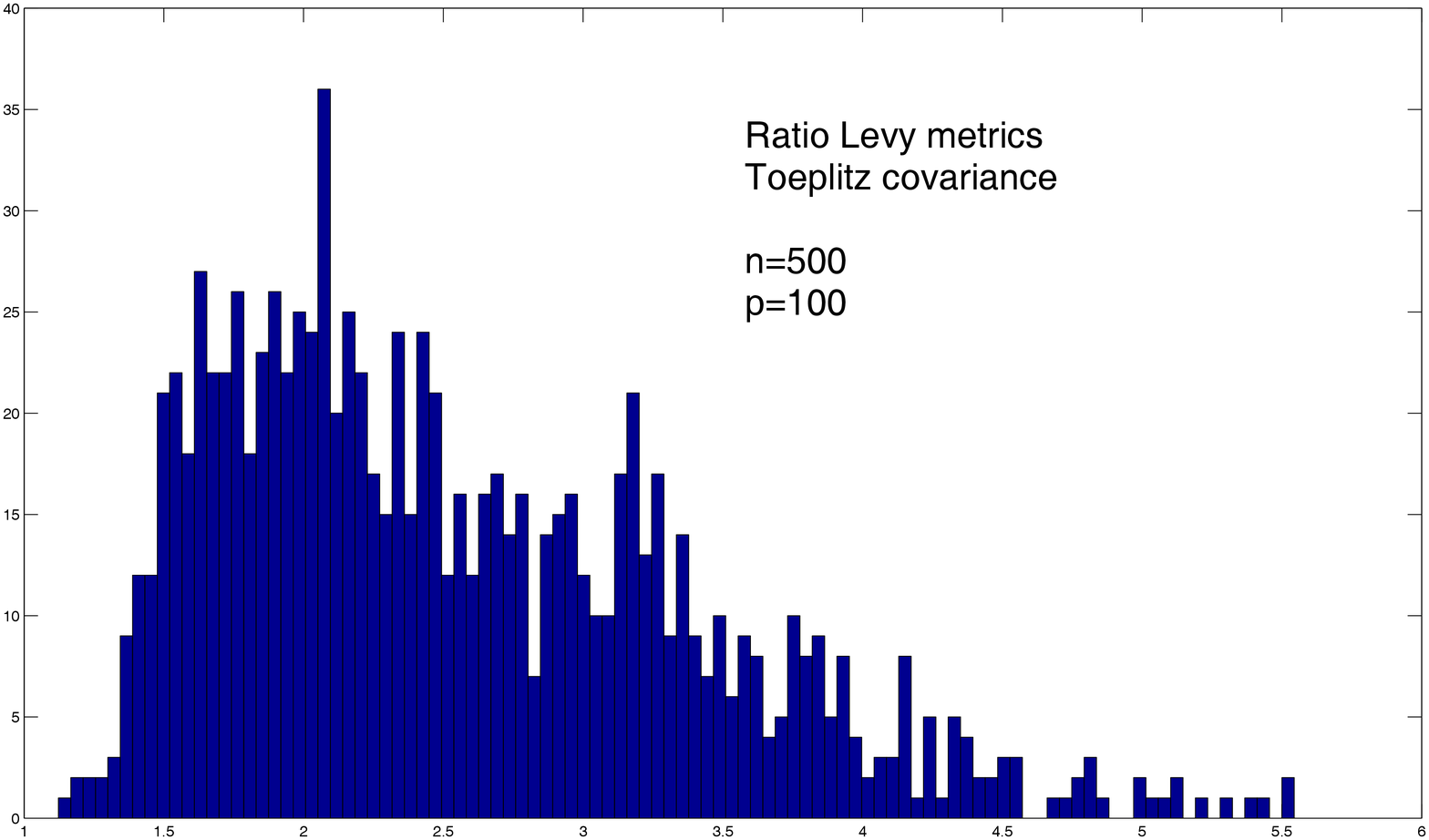}
\end{center}
\caption{\label{fig:ratioLevyToep} \textbf{Case
$\boldsymbol{\Sigma_p}$ Toeplitz with entries ($.3^{|i-j|}$)}:
Ratios $d_L(\hatHp,\Hp)/d_L(\Fp,\Hp)$ over 1,000 repetitions.
Dictionary consisted of only point masses. Large values indicate
better performance of our algorithm. All ratios were found to be
larger than 1.}
\end{figure}

Once again, the results displayed in Figure \ref{fig:Toep500100Case}
are quite encouraging. Note that this time, the population spectral
distribution could only be approximated by a large number of
elements of our dictionary. So there was no sparse representation of
$\Hinfty$ in our chosen dictionary of measures. However, computation
time was not severely affected and the results are still quite good.
To give a more detailed comparison, we present in Figure
\ref{fig:ratioLevyToep} a histogram of ratios $d_L(\hatHp,\Hp)/d_L(\Fp,\Hp)$.

\section{Consistency}\label{sec:Consistency}
In this section, we prove that the algorithm we propose leads to a
consistent (in the sense of weak convergence of probability
measures) estimator of the spectral distribution of the covariance
matrices of interest.

More precisely, we focus on the ``$L_{\infty}$" version of the
algorithm proposed in \ref{subsubsec:convexOptFormulation}. In
short, the theoretical results we prove state that as our
computational resources grow (both in terms of size of available
data and grid points on which to evaluate functions), the estimator
$\hatHp$ converges to $\Hinfty$. The meaning of Theorem
\ref{thm:consistency}, which follows, is the following. We first choose a family of points $\{z_j\}$
in the upper-half of the complex plane, with a limit point in the
upper-half of the complex plane. We assume that the population
spectral distribution $\Hp$ has a limit, in the sense of weak
convergence of distributions, when $p\tendsto \infty$. We call this
limit $\Hinfty$. This assumption of weak convergence allows us to
vary $\Hp$, as $p$ grows, and to not be limited to $\Hp=\Hinfty$ for
the theory; this provides maximal generality. We then solve the
``$L_{\infty}$" version of our optimization problem, by including
more and more of the $z_j$'s in the optimization problem as
$n\tendsto
\infty$. We assume in Theorem
\ref{thm:consistency} that we can solve this problem by optimizing
over all probability measures. Then Theorem
\ref{thm:consistency} shows that the solution of the optimization
problem, $\hatHp$, converges in distribution to the limiting
population spectral distribution, $\Hinfty$. In Corollary
\ref{coro:consistencyDiscCas}, we show that the same conclusion
holds if the optimization is now made over probability measures that
are mixture of point masses, whose locations are on a grid whose
step size goes to 0 with $p$ and $n$. Actually, the requirement is
that the dictionary of measures we use contain these diracs. It can
of course be larger. Hence, Corollary \ref{coro:consistencyDiscCas}
proves consistency of the estimators specifically obtained through
our algorithm. Beside the assumptions of Theorem \ref{thm:mapa}, we
assume that all the spectra of the population covariances are
(uniformly) bounded. That translates into the mild requirement that
the support of all $\Hp$'s be contained in a same compact set. Note
that in the context of asymptotics at fixed spectral distribution,
this is automatically satisfied.

We now turn to a more formal statement of the theorem. The notation
$B(z_0,r)$ denotes the closed ball of center $z_0$ and radius $r$.
Our main theorem is the following.
\begin{theorem}\label{thm:consistency}
Suppose we are under the setup of Theorem \ref{thm:mapa}, $H_p
\Rightarrow H_{\infty}$ and $p/n\tendsto
\gamma$, with $0<\gamma<\infty$. Assume that
the spectra of the $\Sigma_p$'s are uniformly bounded. Let
$J_1,J_2,\ldots, $ be a sequence of integers tending to $\infty$.
Let $z_0\in
\mathbb{C}^+$ and $r\in \mathbb{R}^+$ be such that $B(z_0,r)\subset
\cplus$. Let $z_1,z_2,\ldots$ be a sequence of complex variables
with an accumulation point, all contained in $B(z_0,r)$. Let
$\hatHp$ be the solution of
\begin{equation}\label{eq:SolOptimPblem}
\hatHp=\underset{H}{\argmin} \max_{j\leq J_n}
\left|\frac{1}{\stesdd(z_j)}+z_j-\frac{p}{n} \int \frac{\lambda
dH(\lambda)}{1+\lambda \stesdd(z_j)}\right|\;,
\end{equation}
where $H$ is a probability measure. Then we have
$$
\hatHp \Rightarrow H_{\infty}\;, a.s\;.
$$
\end{theorem}
Before we turn to proving the theorem, we need a few intermediate
results. An important step in the proof is the following analytic
lemma.
\begin{lemma}\label{lemma:SolOptimPbConsistent}
Suppose we have a family $\{z_i\}_{i=1}^{\infty}$ of complex numbers
in $\cplus$, with an accumulation point in $\cplus$. Suppose there
exist a sequence $\{J_i\}_{i=1}^{\infty}$ of integers tending to
$\infty$, a sequence $\{\eps_i\}_{i=1}^{\infty}$ of positive reals
tending to $0$, a sequence $\{p(n)\}_{n=1}^{\infty}$ of integers,
with $p(n)/n\tendsto \gamma \in \mathbb{R}_+^*$, and a sequence of
probability measures $\{\hatHp\}_{p=1}^{\infty}$ such that
\begin{equation}\label{eq:inexactMP}
\forall j\leq J_n\,, \left|\frac{1}{\stesdd(z_j)}+z_j-\frac{p}{n} \int \frac{\lambda
d\hatHp(\lambda)}{1+\lambda \stesdd(z_j)}\right|< \eps_n\;.
\end{equation}
Assume that  $\vinfty$ satisfies
\begin{equation}\label{eq:exactMP}
-\frac{1}{\vinfty(z_j)}=z_j -\gamma \int \frac{\lambda
d\Hinfty(\lambda)}{1+\lambda \vinfty(z_j)}\;,
\end{equation}
for some probability measure $H_{\infty}$.  Assume that
$\stesdd(z_j)\tendsto
\vinfty(z_j)$, and both are analytic in $\cplus$ and from $\cplus$ to $\cplus$.
Further, assume that $|\vinfty(z_j)|<C$ for some $C\in \mathbb{R}$,
and $|\imag{\stesdd(z_j)}|>\delta$, as well as
$|\imag{\vinfty(z_j)}|>\delta$, for some $\delta>0$. Then
$$
\hatHp \Rightarrow \Hinfty\;.
$$
\end{lemma}
\begin{proof}
Since $\vinfty$ satisfies
$$
\frac{1}{\vinfty(z_j)}+z_j -\gamma \int \frac{\lambda
d\Hinfty(\lambda)}{1+\lambda \vinfty(z_j)} = 0\;,
$$
equation (\ref{eq:inexactMP}) reads
$$
\left|\frac{1}{\stesdd(z_j)}-\frac{1}{\vinfty(z_j)}+\left(\gamma-\frac{p}{n}\right)\int\frac{\lambda
d\Hinfty(\lambda)}{1+\lambda \vinfty(z_j)}
+\frac{p}{n}\left(\int\frac{\lambda d\Hinfty(\lambda)}{1+\lambda
\vinfty(z_j)}-\int \frac{\lambda d\hatHp(\lambda)}{1+\lambda
\stesdd(z_j)}\right)\right|<\eps_n\;.
$$
Note that since $|\imag{\stesdd}|>\delta$ and
$|\imag{\vinfty}|>\delta$, and given that
$$\left|\frac{1}{\stesdd}-\frac{1}{\vinfty}\right|\leq
\frac{|\stesdd-\vinfty|}{|\imag{\stesdd}||\imag{\vinfty}|}\;,
$$
we have $|1/\stesdd-1/\vinfty|\tendsto 0$.

Also, because $p/n\tendsto \gamma$, the previous equation implies
that
$$
\int\frac{\lambda d\Hinfty(\lambda)}{1+\lambda
\vinfty(z_j)}-\int \frac{\lambda d\hatHp(\lambda)}{1+\lambda
\stesdd(z_j)}\tendsto 0\;.
$$

Now because $\stesdd(z_j)\tendsto \vinfty(z_j)$, we have
\begin{align*}
\left|\int \frac{\lambda
d\hatHp(\lambda)}{1+\lambda \stesdd(z_j)} -\int \frac{\lambda
d\hatHp(\lambda)}{1+\lambda
\vinfty(z_j)}\right|&= \left|\int \frac{\lambda^2(\vinfty(z_j)-\stesdd(z_j))d\hatHp(\lambda)}
{(1+\lambda
\vinfty(z_j))(1+\lambda \stesdd(z_j))}\right|\\&
\leq
\frac{\left|\stesdd(z_j)-\vinfty(z_j)\right|}{|\imag{\stesdd(z_j)}||\imag{\vinfty(z_j)}|}
\tendsto 0\;.
\end{align*}
So we have
$$
\int \frac{\lambda d\hatHp(\lambda)}{1+\lambda
\vinfty(z_j)}\tendsto \int \frac{\lambda d\Hinfty(\lambda)}{1+\lambda
\vinfty(z_j)}\;.
$$
We remark that for $m\in \cplus$, and $G$ a probability measure on
$\mathbb{R}$, whose Stieltjes transform is denoted by $S_G$,
$$
\int \frac{\lambda dG(\lambda)}{1+\lambda
m}=\frac{1}{m}-\frac{1}{m}\int \frac{dG(\lambda)}{1+\lambda
m}=\frac{1}{m}-\frac{1}{m^2}\int
\frac{dG(\lambda)}{1/m+\lambda}=\frac{1}{m}-\frac{1}{m^2}S_G\left(-\frac{1}{m}\right)\;.
$$
Hence, when the assumptions of the lemma are satisfied, we have
$$
S_{\hatHp}\left(-\frac{1}{\vinfty(z_j)}\right)\tendsto
S_{\Hinfty}\left(-\frac{1}{\vinfty(z_j)}\right)\;.
$$
Now since $\vinfty(z_j)$ satisfies Equation (\ref{eq:exactMP}), we
see that if $\vinfty(z_j)=\vinfty(z_k)$, then $z_j=z_k$. Hence,
$\left\{-1/\vinfty(z_j)\right\}_{j=1}^{\infty}$ is an infinite
sequence of complex numbers in $\cplus$. Moreover, because $\vinfty$
is analytic in $\cplus$, it is continuous, and so
$\left\{-1/\vinfty(z_j)\right\}_{j=1}^{\infty}$ has an accumulation
point. Further, because $|\vinfty(z_j)|<\infty$ and
$\imag{\vinfty(z_j)}>\delta$, this accumulation point is in
$\cplus$.

So under the assumptions of the lemma, we have shown that there
exist an infinite sequence $\{y_j\}_{j=1}^{\infty}$ of complex
numbers in $\cplus$, with an accumulation point in $\cplus$,  such
that
$$
S_{\hatHp}(y_j)\tendsto S_{\Hinfty}(y_j)\;, \forall j \;.
$$
According to \citep{geronimohill03}, Theorem 2, this implies that
$$
\hatHp\Rightarrow \Hinfty\;.
$$
\end{proof}
In the context of spectrum estimation, the intuitive meaning of the
previous lemma is that if for a sequence of complex numbers
$\{z_j\}_{j=1}^{\infty}$ with an accumulation point in $\cplus$, we
can find a sequence of $\hatHp$'s approximately satisfying the \mapa
equation at more and more of the $z_j$'s when $n$ grows, then this
sequence of measures will converge to $\Hinfty$.

We now state and prove a few results that will be needed in the
proof of Theorem \ref{thm:consistency}. The first one is a remark
concerning Stieltjes transforms.
\begin{proposition}\label{prop:LipAndUnifContStT}
The Stieltjes transform, $S_H$, of any probability measure $H$ on
$\mathbb{R}$, is Lipschitz $1/\vmin^2$ on $\cplusabovevmin$.

Hence, if $S_{H_n}(z)\tendsto S_{\Hinfty}(z)$ pointwise, where all
the measures considered are probability measures, the convergence is
uniform on compact subsets of $\cplusabovevmin$.
\end{proposition}
\begin{proof}
We first show the Lipschitz character of $S_H$. We have
$$
S_H(z_1)-S_H(z_2)=\int
\left(\frac{1}{\lambda-z_1}-\frac{1}{\lambda-z_2}\right) dH(\lambda)
=(z_1-z_2)\int \frac{dH(\lambda)}{(\lambda-z_1)(\lambda-z_2)}\;.
$$
Now $|\lambda-z_1|>|\imag{\lambda-z_1}|>\vmin$. So
$$
\left|S_H(z_1)-S_H(z_2)\right|\leq\frac{|z_1-z_2|}{\vmin^2}\;.
$$
So we have shown that $S_H$ is uniformly Lipschitz $1/\vmin^2$ on
$\cplus\bigcap\{\imag{z}>\vmin\}$.

Now, it is an elementary and standard fact of analysis that if a
sequence of $K$-Lipschitz functions converge pointwise to a
$K$-Lipschitz function, then the convergence is uniform on compact
sets. This shows the uniform convergence part of our statement.
\end{proof}
In the proof of the Theorem, we will need the result of the
following proposition.
\begin{proposition}\label{prop:someBounds}
Assume the assumptions underlying Theorem \ref{thm:mapa} are
satisfied. Recall that $\stesdd$ is the Stieltjes transform of
$\tildeFp$, the spectral distribution of $XX^*/n=Y\Sigma_p Y^*/n$.
Assume that the population spectral distribution $H_p$ has a limit
$\Hinfty$ and that all the spectra are uniformly bounded. Let $z \in
B(z_0,r)$, with $B(z_0,r)\subset
\cplus$. Then, almost surely,
$$
\exists N , n>N\Rightarrow \inf_{n,z\in B(z_0,r)} \imag{\stesdd(z)} =
\delta >0\;.
$$
\end{proposition}
\begin{proof}
Since we assume that all spectra are bounded, we can assume that the
population eigenvalues are all uniformly bounded by $K$.  Because
the spectral norm is a matrix norm and $X=Y\Sigma_p^{1/2}$, we have
$$
\lambda_{\mathrm{max}}(X^*X/n)\leq\lambda_{\mathrm{max}}(\Sigma_p)
\lambda_{\mathrm{max}} (Y^*Y/n)\;.
$$
Now it is a standard result in random matrix theory that,
$\lambda_{\mathrm{max}} (Y^*Y/n)\tendsto (1+\gamma)^2$, a.s, so for
$n$ large enough,
$$
\lambda_{\mathrm{max}} (Y^*Y/n)\leq 2 (1+\gamma)^2 \; \text{a.s}\;.
$$
Calling $z=u+iv$, we have
$$
\imag{\stesdd(z)}=\int \frac{vd\tildeFp(\lambda)}{(\lambda-u)^2+v^2}\geq \int
\frac{vd\tildeFp(\lambda)}{2(\lambda^2+u^2)+v^2}\;,
$$
because $v\geq 0$. Now, the remark we made concerning the
eigenvalues of $X^*X/n$ implies that almost surely, for $n$ large
enough, $\tildeFp$ puts all its mass within $[0,C]$, for some $C$.
Therefore,
$$
\imag{\stesdd(z)}\geq \frac{v}{C^2+v^2+2u^2}\;,
$$
and hence $\imag{\stesdd(z)}$ is a.s bounded away from 0, for $n$
large enough.
\end{proof}

To show that we can find ``good" probability measures when solving
our optimization problem, we will need to exhibit a sequence of
measures that approximately satisfy the \mapa equation. The next
proposition is a step in this direction.
\begin{proposition}\label{prop:HinftySolvesProb}
Let $r\in \rplus$ and $z_0 \in \cplus$ be given and satisfying
$B(z_0,r)\subset \cplus$. Suppose $p/n\tendsto \gamma$ when
$n\tendsto \infty$, and $\forall \eps\, \exists N:\, n>N \Rightarrow
\forall z \in B(z_0,r),\; |\stesdd(z)-\vinfty(z)|<\eps$, where $\vinfty$
satisfies equation (\ref{eq:exactMP}). Suppose further that
$|\imag{\vinfty(z)}|>\vmin$ on $B(z_0,r)$. Then, if $\eps<\vmin/2$,
$$
\exists N' \in \mathbb{N}, \; \forall z\in B(z_0,r),\; \forall
n>N', \;
\left|\frac{1}{\stesdd(z)}+z-\frac{p}{n}\int
\frac{\lambda d\Hinfty(\lambda)}{1+\lambda \stesdd(z)}\right|<2\eps
\frac{1+2\gamma}{\vmin^2}
$$
\end{proposition}
\begin{proof}
Using equation (\ref{eq:exactMP}) we find that
\begin{align*}
\Delta_n(z)&=\frac{1}{\stesdd(z)}+z-\frac{p}{n}\int \frac{\lambda
d\Hinfty(\lambda)}{1+\lambda \stesdd(z)}\\
&=\frac{1}{\stesdd(z)}-\frac{1}{\vinfty(z)}+\frac{p}{n}\int
\left(\frac{\lambda}{1+\lambda \vinfty(z)}-\frac{\lambda}{1+\lambda
\stesdd(z)}\right)d\Hinfty(\lambda)\\&+\left(\gamma-\frac{p}{n}\right)\int
\frac{\lambda}{1+\lambda \vinfty(z)}d\Hinfty(\lambda)\\
&\triangleq \Delta_n^{I}(z)+\left(\gamma-\frac{p}{n}\right)\int
\frac{\lambda}{1+\lambda \vinfty(z)}d\Hinfty(\lambda)
\end{align*}
Because $\gamma-p/n\tendsto 0$, and $|\lambda/(1+\lambda
\vinfty(z))|\leq 1/|\imag{\vinfty(z)}|\leq 1/\vmin$, we have
$$
\left(\gamma-\frac{p}{n}\right)\int
\frac{\lambda}{1+\lambda \vinfty(z)}d\Hinfty(\lambda)\tendsto 0 \;\;
\text{uniformly on }B(z_0,r)\;.
$$
Now, of course,
$$
\Delta_n^{I}(z)=\frac{\vinfty(z)-\stesdd(z)}{\stesdd(z)
\vinfty(z)}-\frac{p}{n}(\stesdd(z)-\vinfty(z))\int \frac{\lambda^2}{(1+\lambda
\stesdd(z))(1+\lambda \vinfty(z))} d\Hinfty(\lambda)\;.
$$
We remark that
$|\stesdd(z)|>|\imag{\stesdd(z)}|>\vmin-\eps>\vmin/2$. Hence, if $n$
is large enough,
$$
\left|\Delta_n^{I}(z)\right|\leq 2
\frac{|\vinfty(z)-\stesdd(z)|}{\vmin^2}+2\frac{p}{n}\frac{|\stesdd(z)-\vinfty(z)|}{\vmin^2}\leq
\eps \frac{2}{\vmin^2} \left(1+2\gamma\right)\;.
$$
\end{proof}
We now turn to proving Theorem \ref{thm:consistency}
\begin{proof}[\textbf{Proof of Theorem \ref{thm:consistency}}]
According to Propositions \ref{prop:LipAndUnifContStT} and
\ref{prop:someBounds}, the assumptions put forth in Proposition
\ref{prop:HinftySolvesProb} are a.s satisfied for $\stesdd$ and $\vinfty$ is the Stieltjes
as in Theorem \ref{thm:mapa}. Note also that Theorem \ref{thm:mapa}
states that a.s, $\stesdd(z)\tendsto \vinfty(z)$, and that all these
functions are analytic in $\cplus$. In other words, they have the
properties needed for Lemma \ref{lemma:SolOptimPbConsistent} to
apply.

In particular, Proposition
\ref{prop:HinftySolvesProb} implies that if $\{z_j\}$ is a family of
complex numbers included in $B(z_0,r)$, and if $\hatHp$ is the
solution of equation (\ref{eq:SolOptimPblem}), equation
(\ref{eq:inexactMP}) will be satisfied almost surely, with a family
$\{\eps_j\}$ of positive real numbers that converge to $0$.
According to Lemma
\ref{lemma:SolOptimPbConsistent}, this implies that,
$$
\hatHp\Rightarrow \Hinfty\;, \text{almost surely.}
$$
\end{proof}
As a corollary of Theorem \ref{thm:consistency}, we are now ready to
prove consistency of our algorithm.
\begin{corollary}[Consistency of proposed
algorithm]\label{coro:consistencyDiscCas}
 Assume the same
assumptions as in Theorem \ref{thm:consistency}. Call $\hatHp$ the
solution of equation (\ref{eq:SolOptimPblem}), where the
optimization is now over measures which are sums of atoms, the
location of which are restricted to belong to a grid (depending on
$n$) whose step size is going to $0$ as $n\tendsto \infty$. Then
$$
\hatHp\Rightarrow \Hinfty \; a.s\;.
$$
\end{corollary}
\begin{proof}
All that is needed is to show that a discretized version of
$\Hinfty$ furnishes a good sequence of measures in the sense that
Proposition \ref{prop:HinftySolvesProb} holds for this sequence of
discretized version of $\Hinfty$.

We call $\discHinfty$ a discretization of $\Hinfty$ on a regular
discrete grid of size $1/M_n$. For instance, we can choose
$\discHinfty (x)$ to be a step function, with
$\discHinfty(x)=\Hinfty(x)$ is $x=l/M_n$, $l\in\mathbb{N}$, and
$\discHinfty$ is constant on $[l/M_n,(l+1)/M_n)$. Recall also that
$\Hinfty$ is compactly supported.

In light of the proof of Proposition \ref{prop:HinftySolvesProb},
for the corollary to hold, it is sufficient to show that uniformly
in $z\in B(z_0,r)$,
$$
\left|\int \frac{\lambda}{1+\lambda \stesdd(z)} d\discHinfty(\lambda) - \int \frac{\lambda}{1+\lambda \stesdd(z)}
d\Hinfty(\lambda)\right|\tendsto 0\;.
$$
Now calling $\WassDist{\discHinfty}{\Hinfty}$ the Wasserstein
distance between $\discHinfty$ and $\Hinfty$, we have
$$
\WassDist{\discHinfty}{\Hinfty}=\int_{0}^{\infty}
\left|\discHinfty(x)-\Hinfty(x)\right| dx\tendsto 0 \text{ as } n\tendsto \infty\;.
$$
($\discHinfty$ and $\Hinfty$ put mass only on $\rplus$, so the
previous integral is restricted to $\rplus$. We refer the reader to
the survey \citep{gibbssu01} for properties of different metrics on
probability measures.)

In other respects, it is easy to see that under the assumptions of
Proposition \ref{prop:HinftySolvesProb}, there exists $N$ such that,
$\sup_{n>N, z\in B(z_0,r)} |\stesdd(z)|\leq K$, for some $K<\infty$.
Recall also that under the same assumptions, $\inf_{n>N, z\in
B(z_0,r)}
\imag{\stesdd(z)}\geq \delta$, for some $\delta>0$.

For two probability measures $G$ and $H$, we also have
$$
\WassDist{G}{H}=\sup_f\left\{\left|\int f dG-\int f dH\right|; \; f \text{ a
1-Lipschitz function}\right\}\;.
$$
Hence, because $\Hinfty$ and $\discHinfty$ are supported on a
compact set that is independent of $n$, to have the result we want,
it will be enough to show that
$$
f_{\stesdd(z)}(\lambda)=\frac{\lambda}{1+\lambda \stesdd(z)}
$$
is uniformly Lipschitz (as a function of $\lambda$) when $z\in
B(z_0,r)$ and $n>N$.

Now note that
$$
f_{\stesdd(z)}(\lambda_1)-f_{\stesdd(z)}(\lambda_2)=\frac{\lambda_1-\lambda_2}{(1+\lambda_1
\stesdd(z))(1+\lambda_2 \stesdd(z))}\;.
$$
If $\lambda \leq 1/(2K)$, then $|\lambda \stesdd(z)|\leq 1/2$, so
$|1+\lambda \stesdd(z)|\geq 1/2$. If $\lambda \geq 1/(2K)$, then
$|1+\lambda \stesdd(z)|\geq \lambda \imag{\stesdd(z)}\geq
\delta/(2K)$. So $|1+\lambda \stesdd(z)|\geq
\min(1/2,\delta/(2K))=C$. Hence $f_{\stesdd(z)}$ is
$1/C^2$-Lipschitz, and $C$ is uniform in $n$ and $z$, as needed.

Having thus extended Proposition \ref{prop:HinftySolvesProb} to
discretized versions of $\Hinfty$, the proof of the corollary is the
same as that of Theorem \ref{thm:consistency}.
\end{proof}

The proof of the corollary makes clear that when solving the
optimization problem over any dictionary of probability measures
containing point masses (but also possibly other measures) at grid
points on a grid whose step size goes to $0$, the algorithm will
lead to a consistent estimator.

Finally, as explained in the Appendix, the algorithm we implemented
start with $\stesdd(z_j)$ sequences, as opposed to simply $z_j$
sequences. It can be straightforwardly adapted to handle the $z_j$'s
as a starting point, too, but we got slightly better numerical
results when starting with $\stesdd(z_j)$. The proof we just gave
could be adapted to handle the situation where the $\stesdd(z_j)$'s
are used as  starting point. However, a few other technical issues
would have to be addressed that we felt would make the important
ideas of the proof less clear. Hence we decided to show consistency
in the setting of Corollary
\ref{coro:consistencyDiscCas}.
\section{Conclusion}
In this paper we have presented an original method to estimate the
spectrum of large dimensional covariance matrices. We place
ourselves in a ``large $n$, large $p$" asymptotic framework, where
both the number of observations and the number of variables is going
to infinity, while their ratio goes to a finite, non-zero limit.
Approaching problems in this framework is increasingly relevant as
datasets of larger and larger size become more common.

Instead of estimating individually each eigenvalue, we propose to
associate to each vector of eigenvalues a probability distribution
and estimate this distribution. We then estimate the population
eigenvalues as the appropriate quantiles of the estimated
distribution. We use a fundamental result of random matrix theory,
the \mapa equation, to formulate our estimation problem. We propose
a practical method to solve this estimation problem, using tools
from convex optimization.

The estimator has good practical properties: it is fast to compute
on modern computers (we use the software \citep{mosek} to solve our
optimization problem) and scales well with the number of parameters
to estimate. We show that our estimator of the distribution of
interest is consistent, where the appropriate notion of convergence
is weak convergence of distributions.

The estimator performs a non-linear shrinkage of the sample
eigenvalues. It is basis independent and we hope will help in
improving the estimation of eigenvectors of large dimensional
covariance matrices. To the best of our knowledge, our method is the
first that harnesses deep results of random matrix theory to
practically solve estimation problems. We have seen in simulations
that the improvement it leads to are often dramatic. In particular,
it enables us to find structure in the data when it exists and to
conclude to its absence where there is none, even when classical
methods would point to different conclusions.

\appendix
\begin{center}
    {\bf APPENDIX}
  \end{center}
\setcounter{section}{1}
\subsection{Implementation details}
We plan to release the software we used to create the figures
appearing in the simulation and data analysis section in the near
future. However, we want to mention here the choices of parameters
we made to implement our algorithm. The justifications for them is
based on intuitions coming from studying the equation
(\ref{eq:mapa}).
\paragraph{Scaling of the eigenvalues} If all the entries of the data matrix are
multiplied by a constant $a$, then the eigenvalues of $\Sigma_p$ are
multiplied by $a^2$, and so are the eigenvalues of $\EmpCovMat$.
Hence, if the eigenvalues of $\EmpCovMat$ are divided by a factor
$a$, Equation (\ref{eq:mapa}) remains valid if we change
$H_{\infty}(x)$ into $H_{\infty}(ax)$. In practice, we scale the
empirical eigenvalues by $l_1$ the largest eigenvalue of
$\EmpCovMat$. We solve our convex optimization problem with the
scaled eigenvalues to obtain $H_{\infty}(l_1 x)$, from which we get
$H_{\infty}(x)$ through easy manipulations. The subsequent details
describe how we solve our convex optimization problem, after
rescaling of the eigenvalues.

\paragraph{Choice of $(z_j, v(z_j))$} We have found that using 100 pairs
$(z_j, v(z_j))$ was generally sufficient to obtain good and quick
(10s-60s) results in simulations. More points is of course better.
With 200 points, solving the problem took more time, but was still
doable (40s-3mins). In the simulations and data analysis presented
afterwards, we first chose the $v(z_j)$ and numerically found the
corresponding $z_j$ using \texttt{Matlab}'s optimization toolbox. We
took $v(z_j)$ to have a real part equally spaced (every .02) on
$[0,1]$, and imaginary part of $10^{-2}$ or $10^{-3}$. In other
words, our $v(z_j)$'s consisted of two (discretized) segments in
$\cplus$, the second one being obtained from the first one by a
vertical translation of $9*10^{-3}$.

\paragraph{Choice of interval to focus on} The largest (resp. smallest) eigenvalue of a
$p\times p$ symmetric matrix S are convex (resp. concave) functions
of the entries of the matrix. This is because
$l_1(S)=\sup_{\normhs{u}=1} u'Su$, where $u$ is a vector in
$\mathbb{R}^p$. Hence $l_1(S)$ is the supremum of linear functionals
of the entries of the matrix. Similarly, $l_p(S)=\inf_{\normhs{u}=1}
u'Su$, so $l_p(S)$ is a concave function of the entries of $S$. Note
that the sample covariance matrix $\EmpCovMat$ is an unbiased
estimator of $\Sigma_p$. By Jensen's inequality, we therefore have
$E(l_1(\EmpCovMat))\geq l_1(E(\EmpCovMat))=\lambda_1(\Sigma_p)$. In
other words, $l_1(\EmpCovMat)$ is a biased estimator of
$\lambda_1(\Sigma_p)$, and tends to overestimate it. Similarly,
$l_p(\EmpCovMat)$ is a biased estimator of $\lambda_p(\Sigma_p)$ and
tends to underestimate it. More detailed studies of $l_1$ and $l_p$
indicate that they do not fluctuate too much around their mean.
Practically, as $n\tendsto \infty$, we will have with large
probability, $l_p\leq
\lambda_p$ and $l_1\geq
\lambda_1$. (In certain cases,  concentration bounds can make the previous statement rigorous.)
Hence, after rescaling of the eigenvalues, it will be enough to
focus on probability measures supported on the interval
$[l_p/l_1,1]$ when decomposing $H_{\infty}(l_1 x)$.

\paragraph{Choice of dictionary} In the ``smallest" implementation, we limit ourselves to
a dictionary consisting of point masses on the interval
$[l_p/l_1,1]$, with equal spacing of $.005$. We call $\zeta_p$ the
length of this interval. In larger implementations, we split the
interval $[l_p/l_1,1]$ into dyadic intervals, getting at scale $k$,
$2^k$ intervals:
$[l_p/l_1+j2^{-k}\zeta_p,l_p/l_1+(j+1)2^{-k}\zeta_p]$, for
$j=0,\ldots,2^{k}-1$. We store the end points of all the intervals
at all the scales from $k=2$ to $k=8$ for the coarsest
implementation and up to 10 for the finest We implemented
dictionaries containing:
\begin{enumerate}
\item Point masses every .005 on $[l_p/l_1,1]$, and probability measures supported
on the dyadic intervals described above that have constant density
on these intervals.
\item Point masses every .005 on $[l_p/l_1,1]$, and probability measures supported
on the dyadic intervals described above that have constant density
on these intervals, as well as probability measures on those dyadic
intervals that have linearly increasing and linearly decreasing
densities.
\end{enumerate}
The simulations presented above were made with this latter choice of
dictionary using scales up to 8.

\bibliographystyle{C:/NekTexAuxiliaries/Bibliography/annstats}
\bibliography{C:/NekTexAuxiliaries/Bibliography/research}

\begin{thebibliography}{36}
\expandafter\ifx\csname natexlab\endcsname\relax\def\natexlab#1{#1}\fi
\expandafter\ifx\csname url\endcsname\relax
  \def\url#1{\texttt{#1}}\fi
\expandafter\ifx\csname urlprefix\endcsname\relax\def\urlprefix{URL }\fi

\bibitem[{Akhiezer(1965)}]{akhiezer65}
\textsc{Akhiezer}, N.~I. (1965).
\newblock \emph{The classical moment problem and some related questions in
  analysis}.
\newblock Translated by N. Kemmer. Hafner Publishing Co., New York.

\bibitem[{Anderson(1963)}]{anderson63}
\textsc{Anderson}, T.~W. (1963).
\newblock Asymptotic theory for principal component analysis.
\newblock \emph{Ann. Math. Statist.} \textbf{34}, 122--148.

\bibitem[{Anderson(2003)}]{anderson03}
\textsc{Anderson}, T.~W. (2003).
\newblock \emph{An introduction to multivariate statistical analysis}.
\newblock Wiley Series in Probability and Statistics. Wiley-Interscience [John
  Wiley \& Sons], Hoboken, NJ, third edition.

\bibitem[{Bai(1999)}]{bai99}
\textsc{Bai}, Z.~D. (1999).
\newblock Methodologies in spectral analysis of large-dimensional random
  matrices, a review.
\newblock \emph{Statist. Sinica} \textbf{9}, 611--677.
\newblock With comments by G. J.\ Rodgers and Jack W.\ Silverstein; and a
  rejoinder by the author.

\bibitem[{Baik et~al.(2005)Baik, Ben~Arous, and P\'ech\'e}]{bbap}
\textsc{Baik}, J., \textsc{Ben~Arous}, G., and \textsc{P\'ech\'e}, S. (2005).
\newblock Phase transition of the largest eigenvalue for non-null complex
  sample covariance matrices.
\newblock \emph{Ann. Probab.} \textbf{33}, 1643--1697.

\bibitem[{Baik and Silverstein(2004)}]{baiksilverstein04}
\textsc{Baik}, J. and \textsc{Silverstein}, J. (2004).
\newblock Eigenvalues of large sample covariance matrices of spiked population
  models.
\newblock \emph{arXiv:math.ST/0408165} .

\bibitem[{Bickel and Levina(2004)}]{bickellevina04}
\textsc{Bickel}, P.~J. and \textsc{Levina}, E. (2004).
\newblock Some theory of {F}isher's linear discriminant function, `naive
  {B}ayes', and some alternatives when there are many more variables than
  observations.
\newblock \emph{Bernoulli} \textbf{10}, 989--1010.

\bibitem[{Bickel and Levina(2006)}]{bickellevina06}
\textsc{Bickel}, P.~J. and \textsc{Levina}, E. (2006).
\newblock Regularized estimation of large covariance matrices.
\newblock \emph{Forthcoming Technical Report} .

\bibitem[{B{\"o}ttcher and Silbermann(1999)}]{bottchersilbermann}
\textsc{B{\"o}ttcher}, A. and \textsc{Silbermann}, B. (1999).
\newblock \emph{Introduction to large truncated {T}oeplitz matrices}.
\newblock Universitext. Springer-Verlag, New York.

\bibitem[{Boyd and Vandenberghe(2004)}]{boydvandenberghe04}
\textsc{Boyd}, S. and \textsc{Vandenberghe}, L. (2004).
\newblock \emph{Convex optimization}.
\newblock Cambridge University Press, Cambridge.

\bibitem[{Burda et~al.(2004)Burda, G{\"o}rlich, Jarosz, and
  Jurkiewicz}]{burdaetal04}
\textsc{Burda}, Z., \textsc{G{\"o}rlich}, A., \textsc{Jarosz}, A., and
  \textsc{Jurkiewicz}, J. (2004).
\newblock Signal and noise in correlation matrix.
\newblock \emph{Physica A} \textbf{343}, 295--310.

\bibitem[{Burda et~al.(2005)Burda, Jurkiewicz, and Wac{\l}aw}]{burdajurkwaclaw}
\textsc{Burda}, Z., \textsc{Jurkiewicz}, J., and \textsc{Wac{\l}aw}, B. (2005).
\newblock Spectral moments of correlated {W}ishart matrices.
\newblock \emph{Phys. Rev. E} \textbf{71}.

\bibitem[{Campbell et~al.(1996)Campbell, Lo, and
  MacKinlay}]{CampbellLoMacKinlay}
\textsc{Campbell}, J., \textsc{Lo}, A., and \textsc{MacKinlay}, C. (1996).
\newblock \emph{The Econometrics of Financial Markets}.
\newblock Princeton University Press, Princeton, NJ.

\bibitem[{Chen et~al.(1998)Chen, Donoho, and Saunders}]{chendonohosaunders98}
\textsc{Chen}, S.~S., \textsc{Donoho}, D.~L., and \textsc{Saunders}, M.~A.
  (1998).
\newblock Atomic decomposition by basis pursuit.
\newblock \emph{SIAM J. Sci. Comput.} \textbf{20}, 33--61 (electronic).

\bibitem[{Durrett(1996)}]{durrett96}
\textsc{Durrett}, R. (1996).
\newblock \emph{Probability: theory and examples}.
\newblock Duxbury Press, Belmont, CA, second edition.

\bibitem[{{El Karoui}(To Appear)}]{nekGencov}
\textsc{{El Karoui}}, N. (To Appear).
\newblock Tracy-{W}idom limit for the largest eigenvalue of a large class of
  complex sample covariance matrices.
\newblock \emph{The Annals of Probability} See also \verb=arxiv.PR/0503109=.

\bibitem[{Geman(1980)}]{geman80}
\textsc{Geman}, S. (1980).
\newblock A limit theorem for the norm of random matrices.
\newblock \emph{Ann. Probab.} \textbf{8}, 252--261.

\bibitem[{Geronimo and Hill(2003)}]{geronimohill03}
\textsc{Geronimo}, J.~S. and \textsc{Hill}, T.~P. (2003).
\newblock Necessary and sufficient condition that the limit of {S}tieltjes
  transforms is a {S}tieltjes transform.
\newblock \emph{J. Approx. Theory} \textbf{121}, 54--60.

\bibitem[{Gibbs and Su(2001)}]{gibbssu01}
\textsc{Gibbs}, A.~L. and \textsc{Su}, F. (2001).
\newblock On choosing and bounding probability metrics.
\newblock \emph{International Statistical Review} \textbf{70}, 419--435.

\bibitem[{Gray(2002)}]{gray}
\textsc{Gray}, R.~M. (2002).
\newblock Toeplitz and circulant matrices: A review.
\newblock Available at \verb=http://ee.stanford.edu/~gray/toeplitz.pdf=.

\bibitem[{Grenander and Szeg{\"o}(1958)}]{grenanderszego58}
\textsc{Grenander}, U. and \textsc{Szeg{\"o}}, G. (1958).
\newblock \emph{Toeplitz forms and their applications}.
\newblock California Monographs in Mathematical Sciences. University of
  California Press, Berkeley.

\bibitem[{Hastie et~al.(2001)Hastie, Tibshirani, and Friedman}]{htf01}
\textsc{Hastie}, T., \textsc{Tibshirani}, R., and \textsc{Friedman}, J. (2001).
\newblock \emph{The {E}lements of {S}tatistical {L}earning}.
\newblock Springer Series in Statistics. Springer-Verlag, New York.
\newblock Data mining, inference, and prediction.

\bibitem[{Hiai and Petz(2000)}]{hiaipetz00}
\textsc{Hiai}, F. and \textsc{Petz}, D. (2000).
\newblock \emph{The semicircle law, free random variables and entropy},
  volume~77 of \emph{Mathematical Surveys and Monographs}.
\newblock American Mathematical Society, Providence, RI.

\bibitem[{Johnstone(2001)}]{imj}
\textsc{Johnstone}, I. (2001).
\newblock On the distribution of the largest eigenvalue in principal component
  analysis.
\newblock \emph{Ann. Statist.} \textbf{29}, 295--327.

\bibitem[{Jonsson(1982)}]{jonsson82}
\textsc{Jonsson}, D. (1982).
\newblock Some limit theorems for the eigenvalues of a sample covariance
  matrix.
\newblock \emph{J. Multivariate Anal.} \textbf{12}, 1--38.

\bibitem[{Laloux et~al.(1999)Laloux, Cizeau, Bouchaud, and
  Potters}]{lalouxetal}
\textsc{Laloux}, L., \textsc{Cizeau}, P., \textsc{Bouchaud}, J.-P., and
  \textsc{Potters}, M. (1999).
\newblock Noise dressing of financial correlation matrices.
\newblock \emph{Phys. Rev. Lett.} \textbf{83}, 1467--1470.

\bibitem[{Lax(2002)}]{lax}
\textsc{Lax}, P.~D. (2002).
\newblock \emph{Functional analysis}.
\newblock Pure and Applied Mathematics (New York). Wiley-Interscience [John
  Wiley \& Sons], New York.

\bibitem[{Ledoit and Wolf(2004)}]{LedoitWolf04}
\textsc{Ledoit}, O. and \textsc{Wolf}, M. (2004).
\newblock A well-conditioned estimator for large-dimensional covariance
  matrices.
\newblock \emph{J. Multivariate Anal.} \textbf{88}, 365--411.

\bibitem[{Mar{\v{c}}enko and Pastur(1967)}]{mp67}
\textsc{Mar{\v{c}}enko}, V.~A. and \textsc{Pastur}, L.~A. (1967).
\newblock Distribution of eigenvalues in certain sets of random matrices.
\newblock \emph{Mat. Sb. (N.S.)} \textbf{72 (114)}, 507--536.

\bibitem[{Mardia et~al.(1979)Mardia, Kent, and Bibby}]{mardiakentbibby}
\textsc{Mardia}, K.~V., \textsc{Kent}, J.~T., and \textsc{Bibby}, J.~M. (1979).
\newblock \emph{Multivariate analysis}.
\newblock Academic Press [Harcourt Brace Jovanovich Publishers], London.
\newblock Probability and Mathematical Statistics: A Series of Monographs and
  Textbooks.

\bibitem[{MOSEK(2006)}]{mosek}
\textsc{MOSEK} (2006).
\newblock {MOSEK} {O}ptimization {T}oolbox.
\newblock Available at \verb=www.mosek.com=.

\bibitem[{Paul(To Appear)}]{debashis}
\textsc{Paul}, D. (To Appear).
\newblock Asymptotics of sample eigenstructure for a large dimensional spiked
  covariance model.
\newblock \emph{Statistica Sinica} .

\bibitem[{Silverstein(1995)}]{silverstein95}
\textsc{Silverstein}, J.~W. (1995).
\newblock Strong convergence of the empirical distribution of eigenvalues of
  large-dimensional random matrices.
\newblock \emph{J. Multivariate Anal.} \textbf{55}, 331--339.

\bibitem[{Silverstein and Bai(1995)}]{silversteinbai95}
\textsc{Silverstein}, J.~W. and \textsc{Bai}, Z.~D. (1995).
\newblock On the empirical distribution of eigenvalues of a class of
  large-dimensional random matrices.
\newblock \emph{J. Multivariate Anal.} \textbf{54}, 175--192.

\bibitem[{Wachter(1978)}]{wachter78}
\textsc{Wachter}, K.~W. (1978).
\newblock The strong limits of random matrix spectra for sample matrices of
  independent elements.
\newblock \emph{Ann. Probability} \textbf{6}, 1--18.

\bibitem[{Yin et~al.(1988)Yin, Bai, and Krishnaiah}]{yinbaikrishnaiah88}
\textsc{Yin}, Y.~Q., \textsc{Bai}, Z.~D., and \textsc{Krishnaiah}, P.~R.
  (1988).
\newblock On the limit of the largest eigenvalue of the large-dimensional
  sample covariance matrix.
\newblock \emph{Probab. Theory Related Fields} \textbf{78}, 509--521.

\end{thebibliography}

\end{document}